\documentclass[preprint,12pt]{elsarticle}
\usepackage{amsmath,amssymb,amsthm}
\usepackage{booktabs,longtable,array}
\usepackage[hidelinks]{hyperref}
\usepackage{enumitem}
\newtheorem{theorem}{Theorem}[section]
\newtheorem{lemma}[theorem]{Lemma}
\newtheorem{proposition}[theorem]{Proposition}
\newtheorem{corollary}[theorem]{Corollary}
\theoremstyle{definition}
\newtheorem{definition}[theorem]{Definition}
\newtheorem{example}[theorem]{Example}
\theoremstyle{remark}
\newtheorem{remark}[theorem]{Remark}
\newtheorem{question}[theorem]{Question}

\newcommand{\F}{\mathbb{F}}
\newcommand{\cS}{\mathcal{S}}
\newcommand{\cB}{\mathcal{B}}
\newcommand{\cI}{\mathcal{I}}
\newcommand{\cQ}{\mathcal{Q}}
\newcommand{\Adm}{\mathcal{A}}
\newcommand{\rk}{\operatorname{rank}}
\newcommand{\Span}{\operatorname{span}}
\newcommand{\supp}{\operatorname{supp}}
\newcommand{\GL}{\mathrm{GL}}
\newcommand{\GraphS}[1]{\mathcal{S}(#1)}
\newcommand{\nuo}{\nu_{\circ}}
\newcommand{\bd}{\partial}
\newcommand{\Fz}[1]{\F^{#1}_{0}}
\newcommand{\one}{\mathbf{1}}
\newcommand{\wedgeop}{\mathbin{\wedge}}
\newcommand{\core}{\operatorname{core}}
\newcommand{\diag}{\operatorname{diag}}
\newcommand{\charF}{\operatorname{char}\F}
\newcommand{\Kplus}{K_3^{+}}

\journal{J. Combin. Theory Ser. B}

\begin{document}

\begin{frontmatter}
\title{Representability of Graphical Symplectic Matroids}
\author[iai,rkmveri]{Sayantan Ganguly\corref{cor1}}
\ead{sayantan.ganguly.90@tcgcrest.org}
\cortext[cor1]{Corresponding author.}

\affiliation[iai]{
  organization={
    Institute for Advancing Intelligence (IAI), TCG CREST
  },
  addressline={
Block EM, Plot No. 3, Sector V, Salt Lake
  },
  city={Kolkata},
  postcode={700091},
  state={West Bengal},
  country={India}
}

\affiliation[rkmveri]{
  organization={
    Ramakrishna Mission Vivekananda Educational
    and Research Institute (RKMVERI)
  },
  addressline={Belur},
  city={Howrah},
  postcode={711202},
  state={West Bengal},
  country={India}
}

\begin{abstract}
Chow associated with every finite multigraph $G$ a symplectic matroid $\GraphS{G}$ on the signed edge set $E\sqcup E^{*}$, whose restrictions to transversals are frame matroids of signed graphs, and Sarvepalli asked which of these matroids are representable. For a connected graph containing a cycle, we show that every isotropic representation of $S(G)$ admits an incidence-flow normal form. The remaining representation conditions are expressed through its admissible minors. Consequently, if $\GraphS{G}$ is $\F$-representable, then the bridgeless part of $G$ has a nowhere-zero $(\F,+)$-flow. In characteristic other than two, $\GraphS{G}$ is representable if and only if every non-loop edge on a cycle has a parallel mate. In characteristic two, we characterize the loopless graphs without two vertex-disjoint cycles by the existence of a nowhere-zero $(\F,+)$-flow, and we show that rigidly linked short cycles, and hence subdivisions of the prism, are obstructions. For $3$-connected simple graphs this gives a complete answer: $\GraphS{G}$ is $\F$-representable if and only if $\F$ has characteristic two, $G$ has no two vertex-disjoint cycles, and either $|\F|>2$ or $G\cong K_5$. In particular, the matroid of the Petersen graph is representable over no field. Finally, the admissible corank--nullity polynomial of an arbitrary symplectic matroid is intrinsic and satisfies a three-term contraction identity. For graphical symplectic matroids, this becomes a deletion-contraction recursion that closes up once marked vertices are allowed, and polynomials built from enveloping matroids agree with it exactly for sheltering envelopes.
\end{abstract}


\end{frontmatter}

\section{Introduction}\label{sec:intro}

Symplectic matroids are the Coxeter matroids of type $C_n$ \cite{BorovikGelfandWhite2003}. We define symplectic matroids on the ground set $J=[n]\sqcup[n]^{*}$ with the involution $i\leftrightarrow i^{*}$, and their bases are admissible $k$-subsets of $J$, that is, $k$-subsets containing at most one element of each pair $\{i,i^*\}$. A $k$-dimensional subspace $W\le\F^{2n}$ that is totally isotropic for the standard alternating form, given as the row space of a matrix $H=[X\mid Z]$ with $XZ^{T}=ZX^{T}$, determines a \emph{represented} symplectic matroid: its bases are the admissible $k$-sets $B$ with $\det H[B]\neq 0$. Lagrangian matroids are the rank $n$ symplectic matroids, which are equivalent to Bouchet's delta-matroids. Let $\mathbb{F}_q$ be a finite field. Under the standard symplectic representation, the images of stabilizer groups in $\mathbb{F}_q^{2n}$ are precisely the totally isotropic subspaces; together with a compatible choice of phases, such a subspace determines the corresponding $\mathbb{F}_q$-linear quantum stabilizer code, and this is where the questions answered here come from \cite{Sarvepalli2011}.

Chow \cite{Chow2003} introduced a graphical family of symplectic matroids. For a multigraph $G=(V,E)$ with $E=\{e_1,\dots,e_n\}$, the pair $\{i,i^*\}$ is identified with the edge $e_i$, and $i^*$ is thought of as a \emph{starred} copy of $e_i$. An admissible set $S$ is independent in $\GraphS{G}$ if every component of the spanning subgraph formed by the edges met by $S$ is either a tree or unicyclic with an odd number of starred edges on its cycle. For every transversal $T$, the
restriction $\GraphS{G}|_T$ is the \emph{frame} matroid of the signed graph in which the starred edges are negative; Chow's theorem says that these local frame matroids glue together to a symplectic matroid. In the coding application, a representation of $\GraphS{G}$ over $\F_q$ is a
stabilizer matrix of an $[[\,|E|,\,|E|-r\,]]_q$ quantum code, where $r=r(\GraphS{G})$. Sarvepalli \cite[Sec.~V]{Sarvepalli2011} asked to ``find representations for the graphical symplectic matroids'' or, alternatively, to ``find a criterion to test which of these matroids are
representable''. Chow mentioned that, on the basis of a personal communication with N.~White, he mentioned that, a graph with three vertices and four edges, two of which form a double edge, gives a non-representable example \cite{Chow2003}.

\subsection{Main results}\label{subsec:main-results}
The main difficulty is that a representation of $\GraphS{G}$ is a single isotropic subspace that must realize, all at once, the frame matroids of the $2^{|E|}$ signed graphs obtained by starring edges of $G$. Isotropy is a quadratic condition; in characteristic~$2$, frame matroids of signed graphs are in general not representable canonically. Our starting point is that most of this difficulty disappears once a representation is written in the right coordinates. After that, the problem splits into a flow, a linear system, and determinant conditions.

\begin{enumerate}
\item \emph{A normal form.} Restricted to the unstarred elements, a representation of $\GraphS{G}$ represents the cycle matroid $M(G)$. Since regular matroids are uniquely representable
\cite{BrylawskiLucas1976}, row operations and pair scalings bring every representation to the matrix form \eqref{eq:normalform}, whose first block is the reduced incidence matrix $D$ of $G$ (Theorem~\ref{thm:normalform}). Isotropy gives two linear conditions: $z$ is an $\mathbb F$-flow
and $DZ_1^{\mathsf T}$ is symmetric. The bases further
force $z_e\ne0$ on every cyclic edge; on bridges, the flow condition gives $z_e=0$.
What remains of the matroid is a set of vanishing and nonvanishing minors (Remark~\ref{rem:decision}). The dependent sets with at most two starred elements turn into linear equations, the two-star relations of Lemma~\ref{lem:twostar}.

\item \emph{Flows.} The last row of the normal form is a nowhere-zero $(\F,+)$-flow on the graph
$G^\circ$ obtained by deleting loops and bridges (Theorem~\ref{thm:flow}). For a finite field $\mathbb F_q$, the existence of such a flow is equivalent to the existence of a nowhere-zero integer $q$-flow. So flow theory immediately excludes small fields: $\GraphS{K_4}$ is not
$\F_2$-representable, and the matroid of the Petersen graph is representable neither over $\F_2$ nor over $\F_4$ (Corollary~\ref{cor:flow-small}). Theorem~\ref{thm:oddchar} also excludes $\mathcal{S}(K_4)$, $\mathcal{S}(K_{3,3})$ and
$\mathcal{S}(P_{10})$ in characteristic $\neq2$ see Remark~\ref{rem:flow-char}. 

\item \emph{Characteristic $\neq2$.} We derive a complete answer. The symplectic matroid $\GraphS{G}$ is representable over a field of characteristic $\neq2$ if and only if every edge on a cycle, other than a loop, has a parallel mate (Theorem~\ref{thm:oddchar}); the answer does not depend on the field. For sufficiency we build the representation from signed incidence vectors (Lemma~\ref{lem:signed}), with coefficients on each parallel class that sum to zero. For necessity we contract a forest (Corollary~\ref{cor:loops}(3)) so that an edge without a parallel mate becomes the common edge of a fan of triangles. The two-star relations and isotropy then produce a factor $2$ that must vanish, and three starred edges of a triangle, which form an independent set, end up with proportional columns. In particular, a simple graph with a cycle can only be represented in characteristic~$2$ (Corollary~\ref{cor:simple}).

\item \emph{Characteristic $2$ without two disjoint cycles.} Here the flow obstruction is the only
one: if $G$ is loopless and has no two vertex-disjoint cycles, then $\GraphS{G}$ is
$\F$-representable if and only if $G^\circ$ has a nowhere-zero $(\F,+)$-flow
(Theorem~\ref{thm:char2-nu1}). The representation is the one-row lift \eqref{eq:lift} of such a flow (Proposition~\ref{prop:lift}). The construction works because $G[S]$ cannot have two distinct cyclic components. If its cycle space has dimension at least two, the columns are dependent, and Chow's independence rule also fails.
If that dimension is one, independence is determined by the parity of the starred edges on the unique cycle, so its
independence reduces to a single parity check on that cycle. By Jaeger's theorem \cite{Jaeger1979} the condition always holds when $|\F|\ge8$; over $\F_4$ it asks for a nowhere-zero $4$-flow, and over $\F_2$ for $G^\circ$ to be Eulerian.

\item \emph{Characteristic $2$: obstructions.} Short cycles carry almost no freedom: on a cycle of
length $2$ or $3$ the vectors $w_e$ of \eqref{eq:w} have the rigid form \eqref{eq:rigid}
(Lemma~\ref{lem:rigid}). If two such cycles are linked by other cycles in the sense of
Definition~\ref{def:rigid}, the two rigid forms are incompatible, and $\GraphS{G}$ is not
representable over any field of characteristic $2$ (Theorem~\ref{thm:rigid}). Together with
Theorem~\ref{thm:oddchar}, this shows that the prism, the octahedron and $K_n$ for $n\ge6$ give
matroids representable over no field (Corollary~\ref{cor:prism}). Contracting a forest carries the
obstruction over to long cycles, so every graph that contains a subdivision of the prism is excluded
in characteristic $2$ (Proposition~\ref{prop:prism-subdivision}). The examples in
Example~\ref{ex:char2} show that neither this obstruction nor the flow condition characterizes
characteristic $2$ in general.

\item \emph{$3$-connected graphs.} For $3$-connected graphs the two sides meet. If such a graph has
two vertex-disjoint cycles, Menger's theorem links them by three disjoint paths and
Proposition~\ref{prop:prism-subdivision} applies; otherwise the graph is one of the graphs of Dirac
and Lov\'asz (Remark~\ref{rem:nu1-examples}), for which we write down nowhere-zero flows explicitly
(Lemma~\ref{lem:dirac-flows}). The result is a complete classification
(Theorem~\ref{thm:three-connected}): for a $3$-connected simple graph $G$, $\GraphS{G}$ is
$\F$-representable if and only if $\charF=2$, $G$ has no two vertex-disjoint cycles, and either
$|\F|>2$ or $G\cong K_5$. In particular, the matroid of the Petersen graph is representable over no
field (Corollary~\ref{cor:petersen}).

\item \emph{Small graphs and Chow's remark.} For a connected unicyclic graph with cycle $C$, the
matrix $\left[\begin{smallmatrix}D&D\\0&\one_C\end{smallmatrix}\right]$ represents $\GraphS{G}$ over every field of characteristic $2$, and characteristic $2$ is necessary when $G$ is simple
(Corollary~\ref{cor:unicyclic}). The doubled triangle, which is the graph behind Chow's example of a non-representable graphical symplectic matroid, is representable exactly over the fields of
characteristic $2$ with more than two elements (Corollary~\ref{cor:doubled}). So Chow's remark is correct over $\F_2$ and in characteristic $\neq2$, but not over $\F_4$ (Remark~\ref{rem:chow}).

\item \emph{Rank enumerators.} The admissible corank--nullity polynomial $T_{\cS}$ of \eqref{eq:TS} is defined for every symplectic matroid. For represented matroids the rank of an
admissible set is the rank of its columns (Lemma~\ref{lem:intrinsic}), so $T_{\cS}$, and with it
Sarvepalli's restricted Tutte-Martin polynomial, does not depend on the representation (Proposition~\ref{prop:rep-indep}). It satisfies a three-term contraction identity (Theorem~\ref{thm:threeterm}). For Lagrangian matroids this is the recursion of the $U$-polynomial, with $T_{\cS}(u,v)=U_{\cS}(u,uv)$ \cite{EFLS2024,Larson2025} (Corollary~\ref{cor:lagrangian}). For graphical symplectic matroids, switching at a vertex identifies the two contractions of an edge, and we get $T_{\GraphS{G}}=u^{\delta_e}T_{\GraphS{G\setminus e}}+2\,T_{\GraphS{G/e}}$ for every edge $e$
that is not a loop (Corollary~\ref{cor:graph-recursion}). A loop is different: contracting its
starred element marks its vertex. Allowing marked vertices
(Definition~\ref{def:marked}, Proposition~\ref{prop:marked}) gives a minor-closed class on which the
recursion is complete (Theorem~\ref{thm:marked} and Corollary~\ref{cor:marked-class}). For general symplectic matroids, however, restriction to the remaining pairs need not produce a symplectic matroid (Remark~\ref{rem:deletion-fails}), so the three-term identity does not automatically give a closed recursion within that class. Finally, a polynomial built from an enveloping matroid agrees with $T_{\cS}$ exactly when the envelope is sheltering (Proposition~\ref{prop:envelope}), which every representation provides (Corollary~\ref{cor:sheltering-rep}). Duchamp's delta-matroid has twelve envelopes, none of them sheltering, and they give two different polynomials (Example~\ref{ex:duchamp}).
\end{enumerate}

\subsection{Paper organization} The rest of the paper is organized as follows. In Section~\ref{sec:prelim}, we fix the notations and recall the necessary definitions that will be used throughout the paper. In Section~\ref{sec:tools} we prove the intrinsic-rank lemma, the contraction lemmas and the normal form. In Section~\ref{sec:flow} we discuss the flow obstruction. The classification in characteristic $\neq2$ has been discussed in Section~\ref{sec:odd}, and in Section~\ref{sec:char2} we discuss the results in characteristic~$2$, including the classification for $3$-connected graphs. In Section~\ref{sec:examples} we discuss the unicyclic graphs, the doubled triangle and Chow's remark, and in Section~\ref{sec:enumerators} we discuss the rank enumerators. We propose a few open problems in Section~\ref{sec:open}. Statements whose proofs use a computer search are confined to Examples~\ref{ex:char2} and~\ref{ex:duchamp}. The searches and lists the matrices found, and the programs, are provided as supplementary material in \ref{app:computations}.

\section{Preliminaries}\label{sec:prelim}
We begin by fixing the notation and recalling the necessary definitions.
\subsection{Notation}\label{subsec:notation}

Throughout the paper, $\F$ is a field, $\F^\times=\F\setminus\{0\}$, and $\F_q$ is the field with $q$
elements. For a vector $y$, $\supp(y)$ is the set of coordinates at which $y$ is nonzero; for a
column vector $y$ and a scalar $c$, $(y;c)$ is the column vector obtained from $y$ by appending the
entry $c$. We follow matroid terminologies from Oxley \cite{Oxley2011}, and for symplectic matroids we follow
Borovik, Gelfand and White \cite{BorovikGelfandWhite2003}.

\subsection{Symplectic matroids}\label{subsec:symplectic}

\begin{definition}[Admissible sets]\label{def:admissible}
Let $n\ge0$, $[n]=\{1,\dots,n\}$, $[n]^*=\{1^*,\dots,n^*\}$ and $J=[n]\sqcup[n]^*$, with the
fixed-point-free involution $a\mapsto a^*$ given by $i\mapsto i^*$ and $i^*\mapsto i$. For
$A\subseteq J$ put $A^*=\{a^*:a\in A\}$; for $a\in J$ put $J_a=\{a,a^*\}$, and for $Y\subseteq[n]$
put $J_Y=Y\sqcup Y^*$. A set $A\subseteq J$ is \emph{admissible} if $A\cap A^*=\emptyset$, and a
\emph{transversal} if moreover $|A|=n$. We write $\Adm(J)$ for the family of admissible sets and
$\Adm_k(J)$ for its members of size $k$; for a $*$-closed subset $J'\subseteq J$, $\Adm(J')$ and
$\Adm_k(J')$ are defined in the same way.
\end{definition}

\begin{definition}[Symplectic matroids {\cite[Ch.~3]{BorovikGelfandWhite2003}}]\label{def:symplectic}
A linear order $\prec$ on $J$ is \emph{admissible} if $a\prec b$ implies $b^*\prec a^*$. It induces
the \emph{Gale order} on $k$-subsets of $J$: if $A=\{a_1\succ\dots\succ a_k\}$ and
$B=\{b_1\succ\dots\succ b_k\}$, then $A\le B$ if $a_j\preceq b_j$ for $j=1,\dots,k$. A
\emph{symplectic matroid of rank $k$} on $J$ is a pair $\cS=(J,\cB)$ with
$\emptyset\neq\cB\subseteq\Adm_k(J)$ such that, for every admissible order $\prec$, some
$B_\prec\in\cB$ satisfies $B\le B_\prec$ for all $B\in\cB$. The members of $\cB=\cB(\cS)$ are the
\emph{bases}, the subsets of bases are the \emph{independent sets}, forming the family
$\cI(\cS)$, and $r(\cS)=k$ is the \emph{rank}. An element $a\in J$ is a \emph{loop} of $\cS$ if it
lies in no basis. $\cS$ is \emph{Lagrangian} if $k=n$. The \emph{uniform} symplectic matroid
$U^*_{k,n}$ has $\cB=\Adm_k(J)$.
\end{definition}

The bases of a symplectic matroid are the maximal independent sets, so $\cS$ is determined by
$\cI(\cS)$. Chow characterized the families $\cI(\cS)$.

\begin{theorem}[Chow {\cite{Chow2003}}]\label{thm:chow-axioms}
A nonempty subset-closed family $\cI$ of admissible subsets of $J$ is the family of independent sets
of a symplectic matroid if and only if the following holds. If $I,I'\in\cI$, $|I|<|I'|$ and
$I\cup\{y\}\notin\cI$ for every $y\in I'\setminus I$, then $I\cup I'$ is not admissible, and there is
$x\notin I$ such that $I\cup\{x\}\in\cI$ and $(I\setminus I'^{*})\cup\{x^*\}\in\cI$.
\end{theorem}

\begin{lemma}[Restriction to a transversal]\label{lem:restriction}
Let $\cS$ be a symplectic matroid on $J$ and $T$ a transversal. Then
$\cI(\cS)|_T=\{I\in\cI(\cS):I\subseteq T\}$ is the family of independent sets of an ordinary
matroid $\cS|_T$ on $T$.
\end{lemma}
\begin{proof}
The family is subset-closed and contains $\emptyset$. Let $I,I'\in\cI(\cS)|_T$ with $|I|<|I'|$. If
no $y\in I'\setminus I$ satisfied $I\cup\{y\}\in\cI(\cS)$, then $I\cup I'$ would be inadmissible
by Theorem~\ref{thm:chow-axioms}, which is impossible because $I\cup I'\subseteq T$. Hence the
augmentation axiom holds.
\end{proof}

\begin{definition}[Rank function]\label{def:rank}
The \emph{rank function} of a symplectic matroid $\cS$ is
$r_{\cS}:\Adm(J)\to\mathbb{Z}_{\ge0}$,
\begin{equation}\label{eq:rank}
  r_{\cS}(A):=\max\{|I|:\ I\subseteq A,\ I\in\cI(\cS)\}=\max_{B\in\cB(\cS)}|A\cap B| .
\end{equation}
\end{definition}

The two expressions in \eqref{eq:rank} agree because $A\cap B$ is independent and every
independent $I\subseteq A$ lies in $A\cap B$ for some basis $B$.

\begin{lemma}\label{lem:rank-properties}
Let $\cS$ be a symplectic matroid of rank $k$ and let $A,A'\in\Adm(J)$.
\begin{enumerate}
\item For every transversal $T\supseteq A$, $r_{\cS}(A)$ is the rank of $A$ in the matroid
$\cS|_T$.
\item $0\le r_{\cS}(A)\le\min(k,|A|)$, and $r_{\cS}(A)\le r_{\cS}(A\cup\{a\})\le r_{\cS}(A)+1$
whenever $A\cup\{a\}$ is admissible.
\item If $A\cup A'$ is admissible, then
$r_{\cS}(A\cup A')+r_{\cS}(A\cap A')\le r_{\cS}(A)+r_{\cS}(A')$.
\item $\cI(\cS)=\{A\in\Adm(J):r_{\cS}(A)=|A|\}$ and $\cB(\cS)=\{A\in\Adm_k(J):r_{\cS}(A)=k\}$.
\end{enumerate}
\end{lemma}
\begin{proof}
(1) The independent subsets of $A$ are exactly the independent sets of $\cS|_T$ contained in $A$.
Items (2) and (3) follow from (1) and the rank axioms of the matroid $\cS|_T$, with $T$ a
transversal containing $A\cup\{a\}$, respectively $A\cup A'$. Item (4) is immediate from
\eqref{eq:rank}.
\end{proof}

\begin{remark}\label{rem:rank}
The function $r_{\cS}$ depends on $\cS$ and in particular on its rank $k$. For instance
$r_{U^*_{k,n}}(A)=\min(k,|A|)$. If $k=n$, then $r_{\cS}$ is the function $h_D=(g_D+|\cdot|)/2$ of
\cite[Sec.~2.2]{Larson2025}, where $g_D(A)=\max_B(|A\cap B|-|A\cap B^*|)$ is Larson's bisubmodular
rank function of the delta-matroid $D$ of $\cS$; indeed $|A\cap B|+|A\cap B^*|=|A|$ for every
transversal $B$. In that case $r_{\cS}$ satisfies Bouchet's axiom
$r_{\cS}(A\cup\{i\})+r_{\cS}(A\cup\{i^*\})\ge2r_{\cS}(A)+1$ for $A\cap J_i=\emptyset$
\cite[Prop.~4.2]{Bouchet1997}, \cite[Prop.~2.10]{Larson2025}. For $k<n$ this axiom fails, for
example for $U^*_{k,n}$ and $|A|=k$, since then $r(A\cup\{i\})=r(A\cup\{i^*\})=r(A)=k$. For $k<n$,
\eqref{eq:rank} is simply the rank function of the independence system $\cI(\cS)$; we do not know an
axiomatic characterization of these functions.
\end{remark}

\begin{lemma}[Contraction and truncation]\label{lem:minors}
Let $\cS$ be a symplectic matroid of rank $k\ge1$ on $J$.
\begin{enumerate}
\item Let $a\in J$ be a non-loop of $\cS$. Then $\cS/a:=\big(J\setminus J_a,\ \{B\setminus\{a\}:
a\in B\in\cB(\cS)\}\big)$ is a symplectic matroid of rank $k-1$ with
$\cI(\cS/a)=\{I\subseteq J\setminus J_a:\ I\cup\{a\}\in\cI(\cS)\}$.
\item For $0\le j\le k$, the \emph{truncation} $\cS^{(j)}$, whose bases are the independent sets of
$\cS$ of size $j$, is a symplectic matroid of rank $j$, and $r_{\cS^{(j)}}=\min(r_{\cS},j)$.
\end{enumerate}
\end{lemma}
\begin{proof}
(1) Put $\cI_a=\{I\subseteq J\setminus J_a:\ I\cup\{a\}\in\cI(\cS)\}$. It is subset-closed, consists
of admissible sets, and contains $\emptyset$ because $a$ is not a loop. We verify the condition of
Theorem~\ref{thm:chow-axioms}. Let $I,I'\in\cI_a$ with $|I|<|I'|$ and $I\cup\{y\}\notin\cI_a$ for all
$y\in I'\setminus I$. Put $\hat I=I\cup\{a\}$ and $\hat I'=I'\cup\{a\}$. These are independent in
$\cS$, $|\hat I|<|\hat I'|$, $\hat I'\setminus\hat I=I'\setminus I$, and $\hat I\cup\{y\}\notin\cI(\cS)$
for $y\in\hat I'\setminus\hat I$. By Theorem~\ref{thm:chow-axioms}, $\hat I\cup\hat I'$ is not
admissible; since $a^*\notin I\cup I'$, already $I\cup I'$ is not admissible. Moreover there is
$x\notin\hat I$ with $\hat I\cup\{x\}\in\cI(\cS)$ and $(\hat I\setminus\hat I'^*)\cup\{x^*\}\in\cI(\cS)$.
Here $x\neq a$, and $x\neq a^*$ because $\hat I\cup\{a^*\}$ is not admissible; so $x\in J\setminus J_a$.
Since $a\notin\hat I'^*$, we have $\hat I\setminus\hat I'^*=(I\setminus I'^*)\cup\{a\}$. Hence
$I\cup\{x\}\in\cI_a$ and $(I\setminus I'^*)\cup\{x^*\}\in\cI_a$. By Theorem~\ref{thm:chow-axioms},
$\cI_a$ is the family of independent sets of a symplectic matroid. Its maximal members are the sets
$B\setminus\{a\}$ with $a\in B\in\cB(\cS)$: every $I\in\cI_a$ satisfies $I\cup\{a\}\subseteq B$ for
some basis $B$, which contains $a$, and all these sets have the maximum possible size $k-1$.

(2) Put $\cI_j=\{I\in\cI(\cS):|I|\le j\}$. Let $I,I'\in\cI_j$ with $|I|<|I'|$ and
$I\cup\{y\}\notin\cI_j$ for $y\in I'\setminus I$. Since $|I\cup\{y\}|\le|I'|\le j$, we get
$I\cup\{y\}\notin\cI(\cS)$. Theorem~\ref{thm:chow-axioms} for $\cI(\cS)$ shows that $I\cup I'$ is
not admissible and yields $x$ with $I\cup\{x\}$ and $(I\setminus I'^*)\cup\{x^*\}$ in $\cI(\cS)$.
Both sets have size at most $|I|+1\le j$, so they lie in $\cI_j$. Every member of $\cI_j$ extends
to an independent set of size $j$, because it lies in a basis of $\cS$. The formula for the rank
function follows from \eqref{eq:rank}.
\end{proof}

\subsection{Isotropic representations}\label{subsec:representations}

\begin{definition}[Representations]\label{def:representation}
Let $\F$ be a field. Equip $\F^{J}=\F^{2n}$ with the alternating form
$\omega(x,y)=\sum_{i=1}^{n}(x_iy_{i^*}-x_{i^*}y_i)$. Let $H=[X\mid Z]\in\F^{k\times 2n}$ have rank
$k$, with the columns of $X$ indexed by $[n]$ and those of $Z$ by $[n]^*$. For $A\subseteq J$,
$H[A]$ denotes the submatrix formed by the columns indexed by $A$. We call $H$ \emph{isotropic} if
$XZ^{T}=ZX^{T}$, that is, if its row space $W$ satisfies $\omega(W,W)=0$. An isotropic $H$
\emph{represents} the symplectic matroid $\cS$ if $\cB(\cS)=\{B\in\Adm_k(J):\det H[B]\neq0\}$, and
$\cS$ is \emph{$\F$-representable} if some isotropic matrix over $\F$ represents it.
\end{definition}

For every isotropic $H$ of rank $k$, the family $\{B\in\Adm_k(J):\det H[B]\neq0\}$ is the set of
bases of a symplectic matroid \cite[Ch.~3]{BorovikGelfandWhite2003}; it is nonempty by the proof
of Lemma~\ref{lem:intrinsic} below. For $x,y\in\F^{m}$ we write $x\wedgeop y:=xy^{T}-yx^{T}$. If
$x_i,z_i$ denote the columns of $H$ indexed by $i$ and $i^*$, then
$XZ^T-ZX^T=\sum_{i=1}^{n}x_i\wedgeop z_i$, so $H$ is isotropic if and only if
\begin{equation}\label{eq:isotropy}
  \sum_{i=1}^{n}x_i\wedgeop z_i=0 .
\end{equation}
Left multiplication by $\GL_k(\F)$ and the \emph{pair scalings}
$(x_i,z_i)\mapsto(c^{-1}x_i,cz_i)$ with $c\in\F^\times$ preserve isotropy and the represented
matroid. A permutation of the pairs preserves isotropy and relabels the represented matroid
accordingly. The \emph{signed swap} $(x_i,z_i)\mapsto(-z_i,x_i)$ preserves isotropy and exchanges the
roles of $i$ and $i^*$; in characteristic $2$ it is the plain swap $(x_i,z_i)\mapsto(z_i,x_i)$.

\subsection{Chow's graphical symplectic matroid}\label{subsec:chow}

Throughout, a \emph{graph} is a finite multigraph $G=(V,E)$; loops and parallel edges are allowed
unless stated otherwise. Write $E=\{e_1,\dots,e_n\}$ and identify the pair $\{i,i^*\}$ with the
edge $e_i$. We write $e$ for the unstarred element $i$ and $e^*$ for the \emph{starred} element
$i^*$, so that $J=E\sqcup E^*$. A loop is a cycle of length one, and two parallel edges form a cycle
of length two.

\begin{definition}[Graphical symplectic matroid]\label{def:graphical}
For $S\in\Adm(J)$ let $E[S]$ be the set of edges whose pair meets $S$, and let $G[S]=(V,E[S])$ be
the spanning subgraph they form. An edge $e\in E[S]$ is \emph{starred in $S$} if $e^*\in S$. A
connected component $K$ of $G[S]$ is \emph{good} if $K$ is a tree, or if $|E(K)|=|V(K)|$, that is, $K$ is unicyclic, and the unique cycle $C_K$ contains an odd number of edges starred in $S$. The set $S$ is \emph{independent} if every component of $G[S]$ is good. We write $\cI(G)$ for the family of independent sets and $\GraphS{G}$ for the symplectic matroid with $\cI(\GraphS{G})=\cI(G)$.
\end{definition}

\begin{theorem}[Chow {\cite[Thm.~2]{Chow2003}}]\label{thm:chow-graph}
If $G$ is loopless, then $\cI(G)$ is the family of independent sets of a symplectic matroid.
\end{theorem}

Corollary~\ref{cor:loops} below extends Theorem~\ref{thm:chow-graph} to graphs with loops, so
$\GraphS{G}$ is defined for every graph. Let $\tau(G)$ be the number of components of $G$ that contain no cycle. Here we consider loops as cycles.

\begin{lemma}\label{lem:rankSG}
$r(\GraphS{G})=|V|-\tau(G)$.
\end{lemma}
\begin{proof}
Let $S$ be independent. Each component $K$ of $G[S]$ satisfies $|E(K)|\le|V(K)|$, with
$|E(K)|=|V(K)|-1$ if $K$ is a tree. Each acyclic component of $G$ contains at least one component
of $G[S]$, which is a tree. Summing over the components of $G[S]$ gives
$|S|=|E[S]|\le|V|-\tau(G)$. Conversely, in each component $L$ of $G$ that contains a cycle $C$,
choose an edge $f\in C$ and a spanning tree $Y_L$ of $L$ containing $C\setminus\{f\}$, and put
$S_L=Y_L\cup\{f^*\}$; then $L[S_L]$ is unicyclic with exactly one starred edge on its cycle. In each
acyclic component take a spanning tree. The union of these sets is independent of size
$|V|-\tau(G)$.
\end{proof}

For a transversal $T$, the restriction $\GraphS{G}|_T$ is the frame matroid of the signed graph on
$G$ whose negative edges are the $e$ with $e^*\in T$ \cite{Zaslavsky1982,Zaslavsky1991}, and a
representation of $\GraphS{G}$ has to realize all these frame matroids by a single isotropic
subspace. In the coding application \cite[Sec.~III]{Sarvepalli2011}, an $\F_q$-representation of
$\GraphS{G}$ is a stabilizer matrix of an $[[\,|E|,\,|E|-r(\GraphS{G})\,]]_q$ code.

\begin{lemma}\label{lem:switching}
For a vertex $s$ of $G$, let $\sigma_s:J\to J$ exchange $e$ and $e^*$ for every non-loop edge $e$
incident with $s$, and fix all other elements. Then $S\in\cI(G)$ if and only if
$\sigma_s(S)\in\cI(G)$. Consequently $\sigma_s$ maps $\cB(\GraphS{G})$ onto itself and
$r_{\GraphS{G}}(\sigma_s(A))=r_{\GraphS{G}}(A)$ for all $A\in\Adm(J)$.
\end{lemma}
\begin{proof}
$S$ and $\sigma_s(S)$ meet the same pairs, so $G[S]=G[\sigma_s(S)]$. A cycle of length at least
two contains either zero or exactly two non-loop edges incident with $s$, and a loop is fixed by
$\sigma_s$. Hence the parity of the number of starred edges on every cycle is the same for $S$ and
$\sigma_s(S)$.
\end{proof}

\subsection{Graphs, incidence matrices and flows}\label{subsec:graphs}

For a graph $G$, $c(G)$ is its number of components. A \emph{bridge} is a non-loop edge whose
deletion increases the number of components; an edge is a bridge if and only if it lies on no
cycle. For an edge $e$, $G\setminus e$ is obtained by deleting $e$. If $e$ is not a loop, $G/e$ is
obtained by identifying the ends of $e$ and deleting $e$; edges parallel to $e$ become loops. For
the edge set $P$ of a forest, $G/P$ is obtained by contracting the edges of $P$ one by one. For
$X\subseteq E$ we write $\one_X\in\F^{E}$ for its indicator vector, and for $s\in V$ we write
$\one_s\in\F^{V}$ for the corresponding unit vector. $M(G)$ denotes the cycle matroid of $G$.

\begin{definition}[Incidence matrices]\label{def:incidence}
Fix an orientation of $G$, so that every non-loop edge $e$ has a tail $s$ and a head $t$. The
\emph{oriented incidence matrix} $\widehat D\in\F^{V\times E}$ has column $\one_s-\one_t$ for such an
edge and a zero column for every loop. Fix one vertex, the \emph{root}, in each component of $G$.
The \emph{reduced incidence matrix} $D\in\F^{(|V|-c(G))\times E}$ is obtained from $\widehat D$ by
deleting the rows of the roots. Thus $D$ has exactly one column for each edge $e$; we denote it by
$\bd_e$. Explicitly, $\bd_e$ is $\one_s-\one_t$ with the entries at the roots deleted, and $\bd_e=0$
if $e$ is a loop.
\end{definition}

Reversing the orientation of $e$ replaces $\bd_e$ by $-\bd_e$, and parallel edges with the same
orientation have equal columns. The matrix $D$ has rank $|V|-c(G)$ and represents $M(G)$ over every
field \cite[Sec.~5.1]{Oxley2011}. In the normal form of Section~\ref{sec:tools}, $\bd_e$ is the
upper part of the column of the unstarred element $e$; the starred element $e^*$ contributes a
column of the second block, not of $D$. For $U\subseteq V$ let
\[
  \Fz{U}=\Big\{y\in\F^{V}:\ \supp(y)\subseteq U,\ \sum_{s\in V}y_s=0\Big\}.
\]
If $G$ is connected with root $s_0$, we identify $\F^{V\setminus\{s_0\}}$ with $\Fz{V}$ via
$y\mapsto(y,-\sum_sy_s)$. Then $\bd_e=\one_s-\one_t$, and for every connected subgraph $K$ of $G$
\begin{equation}\label{eq:span}
  \Span\{\bd_e:\ e\in E(K)\}=\Fz{V(K)} .
\end{equation}

\begin{definition}[Flows]\label{def:flows}
Let $\Lambda$ be an abelian group. A \emph{$\Lambda$-flow} on $G$ is a map $z:E\to\Lambda$ such that at
every vertex the sum of the values on the edges leaving it equals the sum of the values on the
edges entering it; loops impose no condition. It is \emph{nowhere-zero} if $z(e)\neq0$ for all $e$.
A \emph{nowhere-zero $k$-flow} is a nowhere-zero $\mathbb{Z}$-flow with values in
$\{\pm1,\dots,\pm(k-1)\}$. For $\Lambda=(\F,+)$, $z$ is a flow if and only if $\widehat Dz=0$, if and only
if $Dz=0$, because the row of a root is minus the sum of the other rows of its component.
\end{definition}

We use the following standard facts. (F1)~Every flow vanishes on bridges. (F2)~For a finite abelian group $\Lambda$, $G$ has a nowhere-zero $\Lambda$-flow if and only if it has a nowhere-zero
$|\Lambda|$-flow \cite{Tutte1954}. (F3)~Every bridgeless graph has a nowhere-zero $\mathbb{Z}_2^3$-flow \cite{Jaeger1979}. (F4)~A loopless graph has a nowhere-zero $\mathbb{Z}_2$-flow if and only if all its degrees are even. (F5)~A cubic graph has a nowhere-zero $4$-flow if and only if it is $3$-edge-colourable, and a nowhere-zero $3$-flow if and only if it is bipartite
\cite{Zhang1997}.

We write $G^{\circ}$ for the graph obtained from $G$ by deleting all loops and all bridges. An edge
of $G^\circ$ lies on a cycle of $G$ that contains no bridge and no loop, so $G^{\circ}$ is
bridgeless and loopless. We write $\nuo(G)$ for the maximum number of pairwise vertex-disjoint
cycles of $G$.

\begin{lemma}[Signed incidence vectors]\label{lem:signed}
Let $\charF\neq2$, and let $K$ be a graph with a sign $\sigma(e)\in\{1,-1\}$ on every edge. For an edge $e$ with ends $s,t$, where $s=t$ if $e$ is a loop, put $c_e=\one_s-\sigma(e)\one_t\in\F^{V(K)}$. The vectors $c_e$, $e\in E(K)$, are linearly independent if and only if every component of $K$ is a
tree or is unicyclic with an odd number of negative edges on its cycle.
\end{lemma}
\begin{proof}
Up to sign, $c_e$ does not depend on the order of $s,t$. Vectors of different components have
disjoint supports, so we may assume that $K$ is connected. If $|E(K)|>|V(K)|$, the vectors are dependent. In a relation $\sum_e\lambda_ec_e=0$, a vertex of degree one forces the coefficient of its edge to vanish, so we may delete leaves repeatedly. This settles the case of trees. If $K$ is unicyclic, we are reduced to its cycle $s_1e_1s_2e_2\cdots s_me_ms_1$ with $m\ge1$, where $e_j$ joins $s_j$ and $s_{j+1}$ (indices modulo $m$; for $m=1$ the cycle is a loop). Write $c_j=\one_{s_j}-\sigma_j\one_{s_{j+1}}$ with $\sigma_j=\sigma(e_j)$. Comparing coefficients at $s_{j+1}$ in $\sum_j\lambda_jc_j=0$ gives $\lambda_{j+1}=\sigma_j\lambda_j$, hence $\lambda_1=(\prod_j\sigma_j)\lambda_1$. If $\prod_j\sigma_j=1$, then
$\lambda_j=\sigma_1\cdots\sigma_{j-1}$ is a nontrivial relation. If $\prod_j\sigma_j=-1$, then $2\lambda_1=0$, so all $\lambda_j$ vanish because $\charF\neq2$.
\end{proof}

\section{The intrinsic rank, contraction, normal form}\label{sec:tools}

\subsection{The rank of admissible sets in a representation}

\begin{lemma}[Intrinsic rank]\label{lem:intrinsic}
Let $H\in\F^{k\times 2n}$ be isotropic of rank $k$, and let $\cS$ be the symplectic matroid it
represents. Then for every $A\in\Adm(J)$
\[ \rk H[A]=r_{\cS}(A). \]
In particular, an admissible set is independent in $\cS$ if and only if its columns are linearly
independent.
\end{lemma}

\begin{proof}
Let $W$ be the row space of $H$, and for $U\subseteq J$ let $L_U\le\F^{J}$ be the span of the unit
vectors indexed by $U$. For a transversal $T$, the projection of $W$ onto the coordinates in $T$
has kernel $W\cap L_{J\setminus T}=W\cap L_{T^*}$, so
\begin{equation}\label{eq:transversal-rank}
 \rk H[T]=k-\dim(W\cap L_{T^*}).
\end{equation}
 
\emph{Step 1.} Let $I$ be admissible with linearly independent columns. Then $W\cap L_{I^*}=0$. Indeed, let $\Omega=\left(\begin{smallmatrix}0&I_n\\-I_n&0\end{smallmatrix}\right)$, so that $\omega(x,y)=x\Omega y^{T}$. For $w\in W$, isotropy gives $H\Omega w^{T}=0$, so $\Omega w^{T}$ is a linear relation among the columns of $H$; its support is $\supp(w)^*$. A nonzero $w\in W\cap L_{I^*}$ would therefore give a nontrivial relation among the columns of $I$.

\emph{Step 2.} Let $U$ be admissible with $W\cap L_U=0$, and let $\{i,i^*\}$ be a pair not met by $U$. Then $W\cap L_{U\cup\{i\}}=0$ or $W\cap L_{U\cup\{i^*\}}=0$. Otherwise choose nonzero $w_1\in W\cap L_{U\cup\{i\}}$ and $w_2\in W\cap L_{U\cup\{i^*\}}$. Since $W\cap L_U=0$, we have $(w_1)_i\neq(w_2)_{i^*}\neq0$, while $(w_1)_{i^*}=(w_2)_i=0$. For $j\neq i$, $U$ contains at most
one of $j,j^*$, so $(w_1)_j(w_2)_{j^*}=(w_1)_{j^*}(w_2)_j=0$. Hence $\omega(w_1,w_2)=(w_1)_i(w_2)_{i^*}\neq0$, contradicting isotropy.

\emph{Step 3.} Let $I\subseteq A$ index a basis of the column space of $H[A]$, so $|I|=\rk H[A]$.
By Step~1, $W\cap L_{I^*}=0$. Applying Step~2 to one pair at a time, extend $I^*$ to a
transversal $T^*$ with $W\cap L_{T^*}=0$. Then $T\supseteq I$, and $\rk H[T]=k$ by
\eqref{eq:transversal-rank}. Extending $I$ inside $T$ to a basis of the column space of $H[T]$ gives $B\in\Adm_k(J)$ with $\det H[B]\neq0$ and $I\subseteq B$. With $A=\emptyset$ this shows that $\cS$ has a basis. Thus $r_{\cS}(A)\ge|A\cap B|\ge|I|=\rk H[A]$. Conversely, the columns of $A\cap B$ are linearly independent for every basis $B$, so $\rk H[A]\ge r_{\cS}(A)$.
\end{proof}

\begin{remark}\label{rem:intrinsic}
For Lagrangian matroids, Bouchet \cite[Thm.~4.1]{Bouchet1988} computed the rank function of a
represented delta-matroid in terms of matrix ranks; see also \cite[Sec.~1]{Larson2025}.
Lemma~\ref{lem:intrinsic} treats isotropic subspaces of every dimension. It also shows that the
ordinary vector matroid of the columns of $H$ assigns to every admissible set the same rank as
$\cS$ (Proposition~\ref{prop:rep-indep}).
\end{remark}

\subsection{Contraction}

\begin{lemma}[Contraction of representations]\label{lem:contraction}
Let $H\in\F^{k\times2n}$ be isotropic and represent $\cS$, and let $a\in J$ have a nonzero column.
After left multiplication by an element of $\GL_k(\F)$, assume that the column of $a$ is the first
unit vector. Let $H/a$ be obtained by deleting the first row and the columns indexed by $a$ and
$a^*$. Then $H/a$ is isotropic of rank $k-1$ and represents $\cS/a$.
\end{lemma}
\begin{proof}
By Lemma~\ref{lem:intrinsic}, $a$ is not a loop, so $\cS/a$ is defined (Lemma~\ref{lem:minors}).
Rows $2,\dots,k$ of $H$ vanish in column $a$. Hence, for $2\le i,j\le k$, the contribution of the
pair $J_a$ to $\omega(h_i,h_j)$ is $h_i(a)h_j(a^*)-h_i(a^*)h_j(a)=0$, and the shortened rows remain
pairwise orthogonal. If a linear combination $w$ of rows $2,\dots,k$ vanished outside $J_a$, then
$w\in W\cap L_{\{a^*\}}$, which is $0$ by Step~1 of the proof of Lemma~\ref{lem:intrinsic} applied
to $I=\{a\}$; hence $H/a$ has rank $k-1$. Expanding along column $a$ gives
$\det H[B'\cup\{a\}]=\pm\det(H/a)[B']$ for every admissible $(k-1)$-set
$B'\subseteq J\setminus J_a$.
\end{proof}

\begin{lemma}[Contraction of graphical symplectic matroids]\label{lem:graph-contraction}
Let $P\subseteq E$ be the edge set of a forest of $G$, viewed as a set of unstarred elements. For
every admissible $S$ disjoint from $J_P$,
\[ S\cup P\in\cI(G)\iff S\in\cI(G/P). \]
\end{lemma}
\begin{proof}
Contracting $P$ maps the components of $G[S\cup P]$ bijectively onto the components of $(G/P)[S]$. A component $K$ containing $p$ edges of $P$ is mapped to a component $K'$ with $|V(K')|=|V(K)|-p$ and $|E(K')|=|E(K)|-p$, because $P\cap E(K)$ is a forest. Hence trees go to trees and unicyclic components to unicyclic components. The cycle of $K'$ consists of the edges of the cycle of $K$ not in $P$, and it has the same starred edges, because the edges of $P$ are unstarred.
\end{proof}

\begin{corollary}\label{cor:loops}
\begin{enumerate}
\item For every graph $G$, $\cI(G)$ is the family of independent sets of a symplectic matroid $\GraphS{G}$.
\item If $P$ is the edge set of a forest of $G$, then $\GraphS{G/P}=\GraphS{G}/P$, the iterated
contraction by the unstarred elements of $P$, and $r(\GraphS{G/P})=r(\GraphS{G})-|P|$.
\item If $H$ represents $\GraphS{G}$, then the iterated contraction $H/P$
(Lemma~\ref{lem:contraction}) represents $\GraphS{G/P}$.
\end{enumerate}
\end{corollary}
\begin{proof}
For loopless $G$, (2) follows from Lemma~\ref{lem:graph-contraction},
Theorem~\ref{thm:chow-graph} and Lemma~\ref{lem:minors}(1), since $P$ is independent and the
independent sets of $\GraphS{G}/P$ are the sets $S$ with $S\cup P\in\cI(G)$. For (1), let $G$ have
loops $\ell_1,\dots,\ell_m$, at vertices $s_1,\dots,s_m$. Let $G'$ be the loopless graph obtained by
replacing each $\ell_j$ by two parallel edges $\ell_j,f_j$ joining $s_j$ to a new vertex $t_j$.
Then $P=\{f_1,\dots,f_m\}$ is a forest and $G'/P=G$, so (1) follows from the loopless case of (2).
Now (2) holds in general by the same argument, and (3) follows from (2) and
Lemma~\ref{lem:contraction}, applied once for each edge of $P$; each such edge has a nonzero column
because $P$ is independent.
\end{proof}

\begin{lemma}[Components]\label{lem:components}
If $G=G_1\sqcup G_2$, then $\GraphS{G}=\GraphS{G_1}\oplus\GraphS{G_2}$, and $\GraphS{G}$ is
$\F$-representable if and only if both $\GraphS{G_1}$ and $\GraphS{G_2}$ are.
\end{lemma}
\begin{proof}
The components of $G[S]$ are those of $G_1[S\cap J_{E(G_1)}]$ and of $G_2[S\cap J_{E(G_2)}]$, which
gives the direct sum decomposition; block-diagonal matrices give one direction of the second claim.
Conversely, let $H$ represent $\GraphS{G}$ and let $B_2$ be a basis of $\GraphS{G_2}$. Then
$\GraphS{G}/B_2$ is $\GraphS{G_1}$ together with the elements of the remaining pairs of $G_2$,
all of which are loops. By Lemmas~\ref{lem:contraction} and~\ref{lem:intrinsic}, the iterated
contraction $H/B_2$ represents $\GraphS{G}/B_2$, and its loop columns are zero. Deleting these pairs
preserves isotropy and leaves a representation of $\GraphS{G_1}$.
\end{proof}

\subsection{The normal form}

\begin{definition}[Normal form]\label{def:normalform}
Let $G$ be a connected graph containing a cycle, with $n=|E|$, a fixed orientation and a fixed root,
and let $D\in\F^{(|V|-1)\times n}$ be its reduced incidence matrix. A representation of $\GraphS{G}$
over $\F$ is \emph{in normal form} if it has the shape
\begin{equation}\label{eq:normalform}
 H=\left[\begin{array}{c|c} D & Z_1\\ 0 & z\end{array}\right]\in\F^{|V|\times 2n},
 \qquad Z_1\in\F^{(|V|-1)\times n},\ z\in\F^{1\times n}.
\end{equation}
We write $Z_{1,e}$ for the column of $Z_1$ indexed by $e$ and $z_e$ for the entry of $z$, so that
the columns of $e$ and $e^*$ are $(\bd_e;0)$ and $(Z_{1,e};z_e)$.
\end{definition}

\begin{theorem}[Normal form]\label{thm:normalform}
Let $G$ be a connected graph containing a cycle. If $\GraphS{G}$ is $\F$-representable, then it has
a representation in normal form; more precisely, every representation of $\GraphS{G}$ over $\F$ is
brought into normal form by left multiplication by an element of $\GL_{|V|}(\F)$ followed by pair
scalings. Every representation \eqref{eq:normalform} in normal form satisfies
\begin{enumerate}
\item[\textup{(N1)}] $Dz^{T}=0$, that is, $z$ is an $\F$-flow on $G$;
\item[\textup{(N2)}] $z_e\neq0$ for every edge $e$ lying on a cycle (loops included), and $z_e=0$
for every bridge $e$;
\item[\textup{(N3)}] $DZ_1^{T}$ is symmetric, that is, $\sum_{e\in E}\bd_e\wedgeop Z_{1,e}=0$.
\end{enumerate}
Conversely, every matrix of the shape \eqref{eq:normalform} with $z\neq0$ that satisfies
\textup{(N1)} and \textup{(N3)} is isotropic of rank $|V|$.
\end{theorem}

\begin{proof}
Let $H'=[X'\mid Z']$ represent $\GraphS{G}$; its rank is $|V|$ by Lemma~\ref{lem:rankSG}. A set
$A\subseteq E$ of unstarred elements has no starred edges, so $A$ is independent in $\GraphS{G}$ if
and only if $G[A]$ is a forest; in particular loops are dependent. By Lemma~\ref{lem:intrinsic},
$X'$ represents the cycle matroid $M(G)$, which has rank $|V|-1$, and loops have zero columns.
Since $M(G)$ is regular, any two of its representations over $\F$ are projectively equivalent \cite{BrylawskiLucas1976}, \cite[Sec.~6.3, 6.6]{Oxley2011}. Applying this to $X'$, viewed as a matrix of rank $|V|-1$ whose column space is identified with $\F^{|V|-1}$, and to $D$, we obtain an injective linear map $\varphi:\F^{|V|-1}\to\F^{|V|}$ and scalars $c_e\in\F^\times$ with $X'_e=c_e\varphi(\bd_e)$ for every non-loop edge $e$. Choose $R\in\GL_{|V|}(\F)$ with $R\varphi(y)=(y;0)$ for all $y$, and apply the pair scalings $(x_e,z_e)\mapsto(c_e^{-1}x_e,c_ez_e)$. These operations preserve isotropy and the represented matroid and produce the matrix form \eqref{eq:normalform}.

For any matrix of this shape,
\[ XZ^T-ZX^T=\begin{pmatrix} DZ_1^{T}-Z_1D^{T} & Dz^{T}\\ -zD^{T} & 0\end{pmatrix}, \]
so isotropy is equivalent to (N1) together with (N3). If $z\neq0$, the last row is not a linear combination of the others, because $D$ has full row rank; this proves the converse.

For (N2), let $e$ lie on a cycle $C$. Choose a spanning tree $Y$ of $G$ containing $C\setminus\{e\}$
(any spanning tree if $e$ is a loop). Then $G[Y\cup\{e^*\}]$ is connected and unicyclic, and its
cycle $C$ contains exactly one starred edge, so $Y\cup\{e^*\}$ is a basis. Expanding along the last
row, $\det H[Y\cup\{e^*\}]=\pm z_e\det D[Y]$, and $D[Y]$ is invertible; hence $z_e\neq0$. Finally,
$z_e=0$ for bridges by (N1) and (F1).
\end{proof}

For an edge $e$ with $z_e\neq0$ we put
\begin{equation}\label{eq:w}
 w_e:=z_e^{-1}Z_{1,e}\in\F^{|V|-1},\qquad\text{so that the column of $e^*$ is } z_e\,(w_e;1).
\end{equation}

Using the next lemma we extract linear information from the matroid.

\begin{lemma}[Two-star relations]\label{lem:twostar}
Let $H$ be a representation in normal form, let $C$ be a cycle of $G$ consisting of non-loop
edges, and let $e\neq f$ be edges of $C$. Then
\[ w_e-w_f\in\Span\{\bd_g:\ g\in C\setminus\{e,f\}\}\ \subseteq\ \Fz{V(C)}.\]
In particular $w_e=w_f$ if $e$ and $f$ are parallel. If $C=\{e,f,g\}$ is a triangle, there is
$t\in\F$ with
\begin{equation}\label{eq:triangle}
 w_e-w_f=t\,\varepsilon_g\bd_g,\qquad w_f-w_g=t\,\varepsilon_e\bd_e,\qquad w_g-w_e=t\,\varepsilon_f\bd_f,
\end{equation}
where $\varepsilon_e\bd_e+\varepsilon_f\bd_f+\varepsilon_g\bd_g=0$, $\varepsilon\in\{\pm1\}$, is the
cycle relation of $C$.
\end{lemma}
\begin{proof}
The set $S=(C\setminus\{e,f\})\cup\{e^*,f^*\}$ is admissible, and the only nontrivial component of
$G[S]$ is the cycle $C$, with two starred edges. So $S$ is dependent. Its unstarred part
$C\setminus\{e,f\}$ is a forest, so the columns $(\bd_g;0)$, $g\in C\setminus\{e,f\}$, are linearly
independent. A linear relation therefore has the form
$\alpha z_e(w_e;1)+\beta z_f(w_f;1)\in\Span\{(\bd_g;0)\}$ with $(\alpha,\beta)\neq(0,0)$. The last
coordinate gives $\alpha z_e=-\beta z_f$, and this is nonzero; dividing by it gives the first
claim. If $C=\{e,f,g\}$ is a triangle, the three relations read $w_e-w_f=\alpha\bd_g$,
$w_f-w_g=\beta\bd_e$ and $w_g-w_e=\gamma\bd_f$. Their sum gives
$\beta\bd_e+\gamma\bd_f+\alpha\bd_g=0$. Since any two of $\bd_e,\bd_f,\bd_g$ are linearly
independent, $(\beta,\gamma,\alpha)$ is a multiple of $(\varepsilon_e,\varepsilon_f,\varepsilon_g)$.
\end{proof}

\begin{remark}\label{rem:decision}
Theorem~\ref{thm:normalform} replaces the $|V|(2n-|V|)$ coordinates of an affine chart of the
Grassmannian by the $(|V|-1)n$ entries of $Z_1$ and a flow $z$. For a finite field $\F$ there are at
most $|\F|^{\,n-|V|+1}$ flows. For a fixed flow, (N3) and the relations obtained as in the proof of
Lemma~\ref{lem:twostar} from all dependent admissible sets with at most two starred elements and a
forest as unstarred part are \emph{linear} in $Z_1$; only the remaining determinant conditions are
nonlinear. This is the procedure used for the computations reported in \ref{app:computations}.
\end{remark}

\section{The flow obstruction}\label{sec:flow}

\begin{theorem}\label{thm:flow}
Let $G$ be a graph and $\F$ a field. If $\GraphS{G}$ is $\F$-representable, then $G^{\circ}$ admits
a nowhere-zero $(\F,+)$-flow. If $\F$ is finite, this is equivalent to a nowhere-zero
$|\F|$-flow on $G^\circ$.
\end{theorem}
\begin{proof}
By Lemma~\ref{lem:components} each component $L$ of $G$ gives a representable $\GraphS{L}$, and
$G^\circ$ is the disjoint union of the graphs $L^\circ$. If $L$ is acyclic, $L^\circ$ has no edges.
Otherwise take a representation of $\GraphS{L}$ in normal form (Theorem~\ref{thm:normalform}). By
(N1) the vector $z$ is a flow on $L$, and by (N2) it is nonzero exactly on the edges lying on
cycles. The edges of $L^\circ$ are the non-loop edges lying on cycles. Deleting loops does not affect
Kirchhoff's law, and deleting bridges does not affect it either, because $z$ vanishes there. So
the restriction of $z$ is a nowhere-zero flow on $L^\circ$. The last statement is (F2).
\end{proof}

\begin{corollary}\label{cor:flow-small}
Let $G$ be loopless.
\begin{enumerate}
\item If $\GraphS{G}$ is $\F_2$-representable, every vertex has even degree in $G^\circ$.
\item If $\GraphS{G}$ is $\F_4$-representable, $G^\circ$ has a nowhere-zero $4$-flow; if moreover
$G^\circ$ is cubic, it is $3$-edge-colourable.
\item $\GraphS{K_4}$ and $\GraphS{K_{3,3}}$ are not $\F_2$-representable, and the matroid
$\GraphS{P_{10}}$ of the Petersen graph is representable neither over $\F_2$ nor over $\F_4$.
\end{enumerate}
\end{corollary}
\begin{proof}
(1) and (2) follow from Theorem~\ref{thm:flow} together with (F2), (F4) and (F5). For (3): $K_4$,
$K_{3,3}$ and $P_{10}$ are cubic and bridgeless, so they have no nowhere-zero $\mathbb{Z}_2$-flow,
and $P_{10}$ is not $3$-edge-colourable.
\end{proof}

\begin{remark}\label{rem:flow-char}
Theorem~\ref{thm:flow} is relevant only in characteristic $2$. Indeed, let $\charF\neq2$ and
suppose that every edge on a cycle of $G$, other than a loop, has a parallel mate, which by
Theorem~\ref{thm:oddchar} below is necessary for representability. Then $G^\circ$ is a union of
parallel classes $\Pi$ with $|\Pi|\ge2$, and choosing values $\lambda_e\in\F^\times$ with
$\sum_{e\in\Pi}\lambda_e=0$ on each class, all edges of a class oriented alike, gives a nowhere-zero
$(\F,+)$-flow on $G^\circ$. Combined with Theorem~\ref{thm:oddchar}, the graphs $K_4$, $K_{3,3}$ and
$P_{10}$ of Corollary~\ref{cor:flow-small} give matroids that are not representable in any characteristic $\neq2$. $\GraphS{K_4}$ and $\GraphS{K_{3,3}}$ are representable over every field of characteristic $2$ with more than two elements (Remark~\ref{rem:nu1-examples}), whereas $\GraphS{P_{10}}$ is representable over no field (Corollary~\ref{cor:petersen}).
\end{remark}

\section{Characteristic different from two: a complete classification}\label{sec:odd}

\begin{theorem}\label{thm:oddchar}
Let $G$ be a graph and $\F$ a field with $\charF\neq2$. Then $\GraphS{G}$ is $\F$-representable if
and only if every edge of $G$ that lies on a cycle and is not a loop has a parallel mate. In
particular, the answer does not depend on $\F$.
\end{theorem}

\begin{proof}
By Lemma~\ref{lem:components} we may assume that $G$ is connected. If $G$ is a tree, every
admissible set is independent, so $\GraphS{G}=U^*_{n,n}$ with $n=|E|$, which is represented by
$[I_n\mid I_n]$. Assume that $G$ contains a cycle. Index rows by $V$, orient every edge that is not
a loop, and give parallel edges the same orientation. For an edge $e$ from $s$ to $t$, where $s=t$
if $e$ is a loop, let $d_e=\one_s-\one_t$ and $q_e=\one_s+\one_t$; thus $d_e=0$ and $q_e=2\one_s$
for a loop. Every edge on a cycle that is not a loop lies in a parallel class $\Pi$ with
$|\Pi|\ge2$, and all edges of such a class lie on cycles. Choose $\lambda_e\in\F^\times$ for $e\in\Pi$
with $\sum_{e\in\Pi}\lambda_e=0$, for instance $(1,-1)$ repeated, together with one triple
$(1,1,-2)$ if $|\Pi|$ is odd; this is possible because $\charF\neq2$. Put $\lambda_e=1$ for every
loop $e$, and
\[ x_e=d_e\ (e\in E),\qquad z_e=\lambda_e q_e\ (e\ \text{not a bridge}),\qquad z_e=d_e\ (e\ \text{a bridge}),\]
and $H=[X\mid Z]\in\F^{V\times 2n}$. Writing $d_\Pi,q_\Pi$ for the common vectors of a class, and
using $d_e=0$ for loops,
\[ \textstyle\sum_e x_e\wedgeop z_e=\sum_{\Pi}\big(\sum_{e\in\Pi}\lambda_e\big)\, d_\Pi\wedgeop q_\Pi+\sum_{e\ \mathrm{bridge}}d_e\wedgeop d_e=0, \]
so $H$ is isotropic by \eqref{eq:isotropy}. Let $S$ be admissible. After rescaling columns, $H[S]$
consists of the vectors $d_e$ for the unstarred edges of $G[S]$ and for its starred edges that
are bridges of $G$, and of the vectors $q_e$ for its starred edges that are not bridges of $G$. These are the vectors $c_e$ of
Lemma~\ref{lem:signed} for the signed graph $G[S]$ whose negative edges are the starred edges that
are not bridges of $G$.
Since bridges of $G$ lie on no cycle, the negative edges on a cycle of $G[S]$ are exactly its
starred edges. By Lemma~\ref{lem:signed}, the columns of $S$ are independent if and only if $S$
satisfies Chow's rule. Finally, the set described in the proof of Lemma~\ref{lem:rankSG} has $|V|$
independent columns, so $\rk H=|V|=r(\GraphS{G})$ and $H$ represents $\GraphS{G}$.

\smallskip
Conversely, suppose that $H$ represents $\GraphS{G}$ over $\F$ with $\charF\neq2$, and that an edge $e_0$ with
distinct ends $a$ and $b$ lies on a cycle and has no parallel mate. By Lemma~\ref{lem:components}
we may assume $G$ is connected.

Let $K_1,\dots,K_m$ be the components of $G-\{a,b\}$. Each of them is adjacent to $a$ or to $b$. Since $e_0$ lies on a cycle and has no parallel mate, at least one $K_j$ is adjacent to both. For each $j$ let $Y_j$ be a spanning tree of $K_j$. Put $P_j=Y_j$ if $K_j$ is adjacent to both $a$ and $b$; otherwise put $P_j=Y_j\cup\{f_j\}$, where $f_j$ is an edge from $K_j$ to the endpoint of $e_0$ to which $K_j$ is adjacent. Then
$P=\bigcup_jP_j$ is a forest. In $G'=G/P$ the vertex set is $\{a,b,k_1,\dots,k_r\}$ with $r\ge1$, where $k_i$ is the image of the $i$-th component adjacent to both $a$ and $b$. The edges of $G'$ are: $e_0$, which is the \emph{only} non-loop edge between $a$ and $b$, because a component merged into $a$ is not adjacent to $b$ and vice versa; for each $i$ a nonempty class $\mathcal{H}_i$ of parallel $a$--$k_i$ edges and a nonempty class $\mathcal{G}_i$ of parallel $k_i$--$b$ edges; and possibly
loops. By Corollary~\ref{cor:loops}(3), $\GraphS{G'}$ is $\F$-representable.

Now orient $e_0$ from $a$ to $b$, the edges of $\mathcal{H}_i$ from $a$ to $k_i$ and those of $\mathcal{G}_i$ from $k_i$ to $b$, choose $b$ as the root, and apply Theorem~\ref{thm:normalform} to $G'$. Every non-loop edge lies on a triangle $ak_ib$, so $z_e\neq0$ and $w_e$ is defined for all of them. Write $\bd^a_i=\one_a-\one_{k_i}$ and $\bd^b_i=\one_{k_i}$ for the common columns of $D$ on $\mathcal{H}_i$ and $\mathcal{G}_i$, with the coordinate of $b$ deleted. Then $\bd_{e_0}=\one_a=\bd^a_i+\bd^b_i$ for every $i$.

By Lemma~\ref{lem:twostar} applied to cycles of length two, $w$
is constant on each $\mathcal{H}_i$ and on each $\mathcal{G}_i$; call the values $w^a_i$ and $w^b_i$.
Apply Lemma~\ref{lem:twostar} to the triangle $\{e_0,h,g\}$ with $h\in\mathcal{H}_i$ and $g\in\mathcal{G}_i$. It gives $w_{e_0}-w^a_i\in\Span(\bd^b_i)$, $w^b_i-w_{e_0}\in\Span(\bd^a_i)$ and $w^a_i-w^b_i\in\Span(\bd^a_i+\bd^b_i)$. Writing $w^a_i=w_{e_0}-t_i\bd^b_i$ and $w^b_i=w_{e_0}+s_i\bd^a_i$, the third relation forces $s_i=t_i$, because $\bd^a_i$ and $\bd^b_i$ are linearly independent. Hence \begin{equation}\label{eq:fan-w} w^a_i=w_{e_0}-t_i\,\bd^b_i,\qquad w^b_i=w_{e_0}+t_i\,\bd^a_i \qquad(t_i\in\F).
\end{equation}

The columns of $D$ corresponding to loops are zero. Put $Z_i=\sum_{h\in\mathcal{H}_i}z_h$ and $Z_i'=\sum_{g\in\mathcal{G}_i}z_g$. Kirchhoff's law at $k_i$ and at $a$ gives $Z_i=Z_i'$ and
$z_{e_0}+\sum_iZ_i=0$. By (N3), \eqref{eq:w} and \eqref{eq:fan-w},
\begin{align*}
0&=\sum_{e}z_e\,\bd_e\wedgeop w_e\\
&= z_{e_0}\,\bd_{e_0}\wedgeop w_{e_0}+\sum_i Z_i\Big[\bd^a_i\wedgeop\big(w_{e_0}-t_i\bd^b_i\big)+\bd^b_i\wedgeop\big(w_{e_0}+t_i\bd^a_i\big)\Big]\\
&= \Big(z_{e_0}+\sum_iZ_i\Big)\bd_{e_0}\wedgeop w_{e_0}-2\sum_i Z_it_i\,\bd^a_i\wedgeop\bd^b_i\\
&= -2\sum_i Z_it_i\,\one_a\wedgeop\one_{k_i}.
\end{align*}
The matrices $\one_a\wedgeop\one_{k_i}$ are linearly independent and $2\neq0$, so $Z_it_i=0$ for
every $i$. Since $z_{e_0}=-\sum_iZ_i\neq0$, some $Z_{i_0}$ is nonzero, and therefore $t_{i_0}=0$.

With $t_{i_0}=0$, \eqref{eq:fan-w} gives $w^a_{i_0}=w^b_{i_0}=w_{e_0}$. Take $h\in\mathcal{H}_{i_0}$ and $g\in\mathcal{G}_{i_0}$. The set $\{e_0^*,h^*,g^*\}$ is a triangle with three starred edges, hence independent in $\GraphS{G'}$. But its three columns are multiples of the same vector $(w_{e_0};1)$, a contradiction.
\end{proof}

\begin{corollary}\label{cor:simple}
If $G$ is a simple graph containing a cycle, then $\GraphS{G}$ is not representable over any field
of characteristic $\neq2$. If every edge of $G$ that lies on a cycle and is not a loop has a
parallel mate, then $\GraphS{G}$ is representable over every field of characteristic $\neq2$.
\end{corollary}

\section{Characteristic two}\label{sec:char2}

Throughout this section $\charF=2$, so signs play no role: $\bd_e=\one_s+\one_t$ for an edge $e$ with
ends $s,t$ (in the coordinates of $\Fz{V}$), and $D$ does not depend on the orientation.

\subsection{The one-row lift and graphs without two disjoint cycles}

\begin{proposition}[One-row lift]\label{prop:lift}
Let $\charF=2$ and let $G$ be a connected loopless graph with $\nuo(G)\le1$ that contains a cycle.
Let $\lambda\in\F^{E}$ be an $\F$-flow on $G$ with $\lambda_e\neq0$ for every edge $e$ that is not a
bridge, and put $d_e=\lambda_e$ if $e$ is not a bridge and $d_e=1$ if $e$ is a bridge. Then
\begin{equation}\label{eq:lift}
 H_\lambda=\left[\begin{array}{c|c} D & D\,\diag(d)\\ 0 & \lambda\end{array}\right]
\end{equation}
represents $\GraphS{G}$ over $\F$.
\end{proposition}
\begin{proof}
$H_\lambda$ has the matrix form \eqref{eq:normalform} with $Z_1=D\diag(d)$ and $z=\lambda\neq0$. (N1) holds because $\lambda$ is a flow, and (N3) holds because $DZ_1^{T}=D\diag(d)D^{T}$ is symmetric. By Theorem~\ref{thm:normalform}, $H_\lambda$ is isotropic of rank $|V|=r(\GraphS{G})$. Let $S$ be admissible. The column of $e^*$ is $(d_e\bd_e;\lambda_e)$; dividing it by $d_e$ turns $H_\lambda[S]$ into $\left(\begin{smallmatrix}D[E[S]]\\ \delta_S^{T}\end{smallmatrix}\right)$, where $\delta_S\in\F^{E[S]}$ is the indicator vector of the edges that are starred in $S$ and are not bridges of $G$. Hence the columns of $S$ are dependent if and only if some nonzero $c\in\F^{E[S]}$
satisfies $D[E[S]]c=0$ and $\delta_S^Tc=0$. The solutions of $D[E[S]]c=0$ are the $\F$-flows on $G[S]$, a space of dimension $\beta=|E[S]|-|V|+c(G[S])$. If $\beta=0$, the columns are independent.
If $\beta=1$, then $G[S]$ has a single cyclic component, which is unicyclic with cycle $C$, and the flows are the multiples of $\one_C$. The edges of $C$ lie on a cycle of $G$ and are not bridges, so $\delta_S^T\one_C$ is the number of starred edges of $C$, taken modulo $2$. Thus the columns are independent if and only if $C$ has an odd number of starred edges. If $\beta\ge2$, the columns are dependent. On the other hand, $G[S]$ cannot have two cyclic components, because $\nuo(G)\le1$. So in each case the columns of $S$ are independent exactly when $S$ satisfies Chow's rule.
\end{proof}

\begin{theorem}\label{thm:char2-nu1}
Let $\charF=2$ and let $G$ be a loopless graph with $\nuo(G)\le1$. Then $\GraphS{G}$ is
$\F$-representable if and only if $G^\circ$ has a nowhere-zero $(\F,+)$-flow. Consequently:
\begin{enumerate}
\item $\GraphS{G}$ is representable over every field of characteristic $2$ with at least $8$
elements;
\item $\GraphS{G}$ is $\F_4$-representable if and only if $G^\circ$ has a nowhere-zero $4$-flow;
\item $\GraphS{G}$ is $\F_2$-representable if and only if every vertex has even degree in $G^\circ$.
\end{enumerate}
\end{theorem}
\begin{proof}
Necessity is Theorem~\ref{thm:flow}. For sufficiency, by Lemma~\ref{lem:components} it suffices to
treat each component $L$ of $G$; note that $\nuo(L)\le1$ and that $L^\circ$ is a union of components
of $G^\circ$. If $L$ is a tree, $\GraphS{L}$ is represented by $[I\mid I]$. Otherwise a nowhere-zero
$\F$-flow on $L^\circ$, extended by $0$ to the bridges of $L$, is a flow $\lambda$ on $L$ as in
Proposition~\ref{prop:lift}, which gives a representation. For (1), $L^\circ$ is bridgeless, so it
has a nowhere-zero $\mathbb{Z}_2^3$-flow by (F3); if $|\F|\ge8$, then $(\F,+)$ is an
$\F_2$-vector space of dimension at least $3$ and contains a copy of $\mathbb{Z}_2^3$. Items (2) and
(3) follow from (F2) and (F4).
\end{proof}

\begin{remark}\label{rem:nu1-examples}
Graphs without two vertex-disjoint cycles were characterized by Lov\'asz \cite{Lovasz1965}. In the $3$-connected case, Dirac \cite{Dirac1963} and Lov\'asz showed that they are $K_5$, $K_5$ minus an edge, the wheels, and the graphs obtained from $K_{3,p}$, $p\ge3$, by adding any set of edges inside
the part of size three. Theorem~\ref{thm:char2-nu1} therefore applies to cycles, theta graphs, flowers of cycles sharing one vertex, wheels, $K_4$, $K_{3,3}$ and $K_5$. For example, $K_5$ is Eulerian, so $\GraphS{K_5}$ is representable over every field of characteristic $2$; $K_4$ and $K_{3,3}$ are cubic and $3$-edge-colourable, so combining Theorem~\ref{thm:oddchar} and Corollary~\ref{cor:flow-small} $\GraphS{K_4}$ and $\GraphS{K_{3,3}}$ are representable exactly over the fields of characteristic $2$ with more than two elements.
\end{remark}

\subsection{An obstruction: rigidly linked short cycles}

\begin{lemma}\label{lem:rigid}
Let $\charF=2$, let $G$ be connected, and let $H$ be a representation of $\GraphS{G}$ in normal
form. For every cycle $C$ of length $2$ or $3$ there are $y_C\in\F^{|V|-1}$ and $t_C\in\F$ such
that
\begin{equation}\label{eq:rigid}
 w_e=y_C+t_C\,\bd_e\qquad\text{for all } e\in C,
\end{equation}
and $t_C=0$ if $|C|=2$.
\end{lemma}
\begin{proof}
If $C=\{e,f\}$, then $w_e=w_f$ by Lemma~\ref{lem:twostar}; take $y_C=w_e$ and $t_C=0$. If
$C=\{e,f,g\}$, then \eqref{eq:triangle} reads $w_e+w_f=t\bd_g$, $w_f+w_g=t\bd_e$,
$w_g+w_e=t\bd_f$, and $\bd_e+\bd_f+\bd_g=0$. Put $y_C=w_e+t\bd_e$ and $t_C=t$. Then
$w_f+t\bd_f=w_e+t(\bd_g+\bd_f)=w_e+t\bd_e=y_C$, and similarly for $g$.
\end{proof}

For cycles of length at least $4$, Lemma~\ref{lem:twostar} gives only
$w_e-w_f\in\Span\{\bd_g:g\in C\setminus\{e,f\}\}$, a space of dimension $|C|-2\ge2$, and no
relation of the form \eqref{eq:rigid} follows.

\begin{definition}[Rigidly linked cycles]\label{def:rigid}
Let $C_1,C_2$ be vertex-disjoint cycles of $G$. The \emph{linking family} $\cQ(C_1,C_2)$ is the set
of cycles $Q$ of $G$ with $E(Q)\cap E(C_1)\neq\emptyset$ and $E(Q)\cap E(C_2)\neq\emptyset$. If
$\cQ(C_1,C_2)\neq\emptyset$, the \emph{core} is $\core(C_1,C_2)=\bigcap_{Q\in\cQ(C_1,C_2)}V(Q)$. The
pair $C_1,C_2$ is \emph{rigidly linked} if both cycles have length $2$ or $3$,
$\cQ(C_1,C_2)\neq\emptyset$, and $\core(C_1,C_2)\subseteq V(C_1)$ or
$\core(C_1,C_2)\subseteq V(C_2)$; for instance, if $\core(C_1,C_2)=\emptyset$.
\end{definition}

\begin{theorem}\label{thm:rigid}
If $G$ contains a rigidly linked pair of cycles, then $\GraphS{G}$ is not representable over any
field of characteristic $2$.
\end{theorem}
\begin{proof}
Let $C_1,C_2$ be rigidly linked, say $\core(C_1,C_2)\subseteq V(C_1)$. They lie in one component of
$G$, so by Lemma~\ref{lem:components} we may assume $G$ connected. Suppose that $\GraphS{G}$ is
$\F$-representable with $\charF=2$, and take a representation in normal form. Every edge used below
lies on a cycle, so $z_e\neq0$ and $w_e$ is defined. We use the coordinates of $\Fz{V}$, so that
\eqref{eq:span} holds. Write $y_j=y_{C_j}$ and $t_j=t_{C_j}$ as in Lemma~\ref{lem:rigid}, and put
$y=y_1+y_2$.

Let $Q\in\cQ(C_1,C_2)$, and pick $h_1\in E(Q)\cap E(C_1)$ and
$h_2\in E(Q)\cap E(C_2)$. The cycle $Q$ has length at least $2$, so it consists of non-loop edges, and Lemma~\ref{lem:twostar} gives $w_{h_1}+w_{h_2}\in\Fz{V(Q)}$. Since $\bd_{h_1},\bd_{h_2}\in\Fz{V(Q)}$, \eqref{eq:rigid} gives
$y=w_{h_1}+w_{h_2}+t_1\bd_{h_1}+t_2\bd_{h_2}\in\Fz{V(Q)}$. Hence
$y\in\bigcap_{Q}\Fz{V(Q)}=\Fz{\core(C_1,C_2)}\subseteq\Fz{V(C_1)}$.

Fix $h_j\in E(C_j)$ and consider
\[ S=\{h_1^*,h_2^*\}\cup(C_1\setminus\{h_1\})\cup(C_2\setminus\{h_2\}). \]
The nontrivial components of $G[S]$ are $C_1$ and $C_2$, each a cycle with exactly one starred edge, so $S$ is independent. On the other hand,
\[ z_{h_1}^{-1}\,\mathrm{col}(h_1^*)+z_{h_2}^{-1}\,\mathrm{col}(h_2^*)=\big(w_{h_1}+w_{h_2};\,0\big)
 =\big(y+t_1\bd_{h_1}+t_2\bd_{h_2};\,0\big). \]
By \eqref{eq:span}, $y$ and $\bd_{h_1}$ lie in $\Fz{V(C_1)}=\Span\{\bd_g:g\in C_1\setminus\{h_1\}\}$,
and $\bd_{h_2}\in\Span\{\bd_g:g\in C_2\setminus\{h_2\}\}$. This is a nontrivial linear relation among the columns of $S$, a contradiction.
\end{proof}

\begin{corollary}\label{cor:prism}
Suppose that $G$ contains vertex-disjoint triangles $x_0x_1x_2$ and $y_0y_1y_2$ and pairwise vertex-disjoint paths $P_0,P_1,P_2$, where $P_i$ joins $x_i$ to $y_i$ and has no inner vertex on the triangles. Then $\GraphS{G}$ is not representable over any field of characteristic $2$. If, in addition, some edge of $G$ that lies on a cycle and is not a loop has no parallel mate (for instance, if $G$ is simple), then $\GraphS{G}$ is representable over no field. This applies to the triangular prism $K_3\square K_2$, the octahedron $K_{2,2,2}$ and $K_n$ for $n\ge6$.
\end{corollary}
\begin{proof}
For $\{i,j,l\}=\{0,1,2\}$ let $Q_{ij}$ be the cycle formed by the edge $x_ix_j$, the path $P_j$,
the edge $y_jy_i$ and the path $P_i$. Then $Q_{01},Q_{12},Q_{20}\in\cQ$. The vertices $x_i$, $y_i$ and the inner
vertices of $P_i$ lie only on the two cycles $Q_{ij}$ and $Q_{li}$ among these three, so no vertex
lies on all of them. Hence the core is empty, the triangles are rigidly linked, and
Theorem~\ref{thm:rigid} applies. The second statement follows from
Theorem~\ref{thm:oddchar}. The octahedron contains two disjoint triangles joined by a perfect
matching, and $K_n$ contains the octahedron for $n\ge6$.
\end{proof}

\begin{example}\label{ex:char2}
The following examples show that neither Theorem~\ref{thm:char2-nu1} nor Theorem~\ref{thm:rigid} characterizes representability in characteristic $2$ when $\nuo(G)\ge2$. Vertices are
$0,1,2,\dots$, and $st$ denotes an edge with ends $s,t$. Items (2), (3) and (5) rely on computer searches, described in \ref{app:computations}, where the matrices found are listed.
\begin{enumerate}
\item \emph{The dumbbell} $X_1$ consists of the triangles $T_1=\{01,12,20\}$ and $T_2=\{34,45,53\}$
joined by the bridge $g=03$. Here $\nuo(X_1)=2$, and $X_1^\circ$ is Eulerian. $\GraphS{X_1}$ is
$\F_2$-representable: take the root $0$, $z=\one_{T_1\cup T_2}$, $Z_{1,e}=\bd_e$ for
$e\in T_1\cup\{g\}$ and $Z_{1,e}=\bd_e+\bd_g$ for $e\in T_2$. With the columns ordered
$01,12,20,34,45,53,03$ and the rows indexed by $1,\dots,5$ and the extra row,
\[
H=\left[\begin{array}{ccccccc|ccccccc}
1&1&0&0&0&0&0& 1&1&0&0&0&0&0\\
0&1&1&0&0&0&0& 0&1&1&0&0&0&0\\
0&0&0&1&0&1&1& 0&0&0&0&1&0&1\\
0&0&0&1&1&0&0& 0&0&0&1&1&0&0\\
0&0&0&0&1&1&0& 0&0&0&0&1&1&0\\
0&0&0&0&0&0&0& 1&1&1&1&1&1&0
\end{array}\right].
\]
(N1) holds because every vertex meets $0$ or $2$ edges of $T_1\cup T_2$, and (N3) holds because
$\sum_{e\in T_2}\bd_e\wedgeop(\bd_e+\bd_g)=\big(\sum_{e\in T_2}\bd_e\big)\wedgeop\bd_g=0$. To see that
$H$ represents $\GraphS{X_1}$, let $S$ be admissible and let $c\in\F_2^{E[S]}$ be the coefficient
vector of a linear relation among the columns of $S$; put $c_g=0$ if $g\notin E[S]$. The columns of
$g$ and $g^*$ both equal $(\bd_g;0)$, and the column of $e^*$ is $(\bd_e;1)$ for $e\in T_1$ and
$(\bd_e+\bd_g;1)$ for $e\in T_2$. Let $\gamma$ be the sum of the $c_e$ over the edges $e\in T_2$
that are starred in $S$, and let $k_j$ be the number of edges of $T_j$ starred in $S$. The first
five rows say that $c+\gamma\one_g$ is a flow on $X_1[S]\cup\{g\}$. As $g$ is a bridge, this flow
vanishes at $g$, so $c_g=\gamma$, and on $T_1\cup T_2$ it equals $\mu_1\one_{T_1}+\mu_2\one_{T_2}$,
where $\mu_j=0$ unless $T_j\subseteq E[S]$. Hence $\gamma=\mu_2k_2$, which must vanish if
$g\notin E[S]$, and the last row reads $\mu_1k_1+\mu_2k_2=0$. Since $c$ is determined by
$(\mu_1,\mu_2)$, a nontrivial relation exists if and only if some $T_j\subseteq E[S]$ has an even
number of starred edges, or $g\in E[S]$ and $T_1\cup T_2\subseteq E[S]$. These are exactly the
admissible sets that violate Chow's rule. By contrast, for every one-row lift \eqref{eq:lift}, two
vertex-disjoint triangles, each with one starred edge, form a dependent set (the case $\beta=2$ in
the proof of Proposition~\ref{prop:lift}); shifting $Z_1$ by $\bd_g$ on one side of the bridge
removes this dependence.
\item \emph{Two triangles joined by two disjoint paths}: $X_2$ has edges
$01$, $12$, $20$, $34$, $45$, $53$, $03$, $25$. Here $\core(T_1,T_2)=\{0,2,3,5\}$ meets both triangles, so
Theorem~\ref{thm:rigid} does not apply, and $\GraphS{X_2}$ is representable over $\F_4$ and over
$\F_8$. It is not $\F_2$-representable, since vertex $0$ has degree $3$ in $X_2=X_2^\circ$.
\item \emph{A $4$-cycle with two opposite edges doubled}: $X_3$ has edges $01$, $01$, $23$, $23$, $02$, $13$.
The two cycles of length two are vertex-disjoint and their core is $\{0,1,2,3\}$.
$\GraphS{X_3}$ is $\F_4$-representable.
\item \emph{A triangle and a cycle of length two, rigidly linked}: $X_4$ has edges $01$, $12$, $20$, two
parallel edges $34$, and $03$, $14$, $23$, $04$. The core of the triangle and the $2$-cycle is $\{3,4\}$,
so they are rigidly linked, and $\GraphS{X_4}$ is not representable over any field of
characteristic $2$; since the triangle edges have no parallel mates, it is representable over no
field. Nevertheless $X_4$ is bridgeless and has a nowhere-zero $4$-flow, so
Theorem~\ref{thm:flow} does not exclude $\F_4$.
\item \emph{A chain of three triangles}: $X_5$ has edges $01$, $12$, $20$, $23$, $34$, $42$, $45$, $56$, $64$. It is
bridgeless and Eulerian, every block has $\nuo\le1$, and no cycle meets both $T_1=\{01,12,20\}$ and
$T_3=\{45,56,64\}$. $\GraphS{X_5}$ is $\F_4$-representable but \emph{not} $\F_2$-representable.
Thus over $\F_2$ the flow condition is not sufficient even when every block has $\nuo\le1$. Since
the subgraph formed by the first two triangles has $\nuo=1$ and is Eulerian, its matroid is
$\F_2$-representable by Theorem~\ref{thm:char2-nu1}, as is that of a triangle; so
$\F_2$-representability is not preserved under $1$-sums.
\end{enumerate}
\end{example}

\subsection{Prism subdivisions and $3$-connected graphs}\label{subsec:three-connected}

Contracting a forest shortens cycles, and Corollary~\ref{cor:loops}(3) carries representations of
$\GraphS{G}$ over to $\GraphS{G/F}$. This turns Corollary~\ref{cor:prism} into an obstruction
involving long cycles.

\begin{proposition}\label{prop:prism-subdivision}
Let $C_1$ and $C_2$ be vertex-disjoint cycles of $G$ of length at least $3$, and let $P_0,P_1,P_2$ be
pairwise vertex-disjoint paths such that $P_i$ joins a vertex $x_i$ of $C_1$ to a vertex $y_i$ of
$C_2$ and has no inner vertex on $C_1\cup C_2$. Then $\GraphS{G}$ is not representable over any field
of characteristic $2$. In particular, this holds whenever $G$ contains a subdivision of the prism
$K_3\square K_2$.
\end{proposition}
\begin{proof}
The vertices $x_0,x_1,x_2$ divide $C_1$ into three paths. Keep one edge on each of them, and let
$F_1$ be the set of the other edges of $C_1$. Define $F_2\subseteq E(C_2)$ in the same way with
respect to $y_0,y_1,y_2$, and let $F_3$ consist of the edges of $P_0\cup P_1\cup P_2$ other than the
edge of $P_i$ at $y_i$, for $i=0,1,2$. Since the paths are disjoint and their inner vertices avoid
$C_1\cup C_2$, the set $F=F_1\cup F_2\cup F_3$ is the edge set of a forest, and $F$ joins every vertex
of $C_1\cup C_2\cup P_0\cup P_1\cup P_2$ to exactly one of the six vertices $x_i,y_i$. In $G/F$ the
images of $C_1$ and $C_2$ are therefore vertex-disjoint triangles, and the edges of the paths $P_i$
at $y_i$ form a perfect matching between them; the other edges of $G$ remain, possibly as loops or
parallel edges. By Corollary~\ref{cor:prism}, $\GraphS{G/F}$ is not representable over a field of
characteristic $2$, and by Corollary~\ref{cor:loops}(3) neither is $\GraphS{G}$. A subdivision of the
prism consists of two cycles and three paths of this kind.
\end{proof}

For $3$-connected graphs, Proposition~\ref{prop:prism-subdivision} and
Theorem~\ref{thm:char2-nu1} together settle the question. The flows needed in
Theorem~\ref{thm:char2-nu1} can be written down explicitly.

\begin{lemma}\label{lem:dirac-flows}
Let $\charF=2$ and $|\F|>2$. Every $3$-connected simple graph without two vertex-disjoint cycles has
a nowhere-zero $(\F,+)$-flow.
\end{lemma}
\begin{proof}
Choose distinct $a,b\in\F^\times$ and put $c=a+b$. Then $a,b,c$ are nonzero and pairwise distinct,
and $a+b+c=0$. In characteristic $2$ an $(\F,+)$-flow is an assignment of values to the edges whose
sum at every vertex is $0$. By Remark~\ref{rem:nu1-examples}, the graph is $K_5$, a wheel, or a graph
$K_{3,p}+L$ obtained from $K_{3,p}$, $p\ge2$, by adding a set $L$ of edges inside the part
$\{s_1,s_2,s_3\}$; the case $p=2$, $|L|=3$ is $K_5$ minus an edge. On $K_5$ the constant value $a$
is a flow, because all degrees are $4$.

Let the wheel have hub $h$ and rim $v_1v_2\cdots v_mv_1$, $m\ge3$. Let $x_1,\dots,x_m$ alternate
between $a$ and $b$, and replace $x_m$ by $c$ if $m$ is odd; then $x_i\neq x_{i+1}$ for all $i$, with
indices modulo $m$. Give the rim edge $v_iv_{i+1}$ the value $x_i$ and the spoke $hv_i$ the value
$x_{i-1}+x_i\neq0$. The sum at $v_i$ is $x_{i-1}+x_i+(x_{i-1}+x_i)=0$, and the sum at $h$ is
$\sum_i(x_{i-1}+x_i)=0$.

On $K_{3,p}+L$, with second part $\{t_1,\dots,t_p\}$, give every edge of $L$ the value $a$, and put
$\gamma_i=a$ if $s_i$ has odd degree in $L$ and $\gamma_i=0$ otherwise. It remains to find a
$3\times p$ matrix with entries in $\{a,b,c\}$, whose $(i,j)$ entry is the value of $s_it_j$, such that
every column sums to $0$ and the $i$th row sums to $\gamma_i$. Since $L$ has an even number of
vertices of odd degree, we may renumber so that $(\gamma_1,\gamma_2,\gamma_3)$ is $(0,0,0)$ or
$(a,a,0)$. Start with
\[
\begin{pmatrix}a&a\\ b&b\\ c&c\end{pmatrix},\qquad
\begin{pmatrix}a&b&c\\ b&c&a\\ c&a&b\end{pmatrix},\qquad
\begin{pmatrix}b&c\\ c&b\\ a&a\end{pmatrix},\qquad
\begin{pmatrix}b&a&b\\ c&c&a\\ a&b&c\end{pmatrix}
\]
if the target is $(0,0,0)$ and $p$ is even or odd, or $(a,a,0)$ and $p$ is even or odd,
respectively, and append pairs of equal columns $(a,b,c)^{T}$ until there are $p$ columns. Every
column sums to $0$, and the row sums are $\gamma_1,\gamma_2,\gamma_3$.
\end{proof}

\begin{theorem}\label{thm:three-connected}
Let $G$ be a $3$-connected simple graph and $\F$ a field. Then $\GraphS{G}$ is $\F$-representable if
and only if $\charF=2$, $G$ has no two vertex-disjoint cycles, and either $|\F|>2$ or $G\cong K_5$.
\end{theorem}
\begin{proof}
A $3$-connected simple graph contains a cycle, and none of its edges has a parallel mate, so
Theorem~\ref{thm:oddchar} excludes every characteristic $\neq2$. Let $\charF=2$, and suppose first
that $G$ has vertex-disjoint cycles $C_1$ and $C_2$; both have length at least $3$. If a set $Z$ of at
most two vertices met every path from $V(C_1)$ to $V(C_2)$, then $G-Z$ would contain a vertex of each
cycle, because $|V(C_j)|\ge3$, but no path between them, contradicting $3$-connectivity. By Menger's
theorem \cite[Thm.~3.3.1]{Diestel2017} there are three disjoint $V(C_1)$--$V(C_2)$ paths, that is,
paths meeting $C_1$ only in their first and $C_2$ only in their last vertex, and
Proposition~\ref{prop:prism-subdivision} shows that $\GraphS{G}$ is not $\F$-representable.

Suppose now that $G$ has no two vertex-disjoint cycles. Since $G$ has no loops and no bridges,
$G^\circ=G$, and by Theorem~\ref{thm:char2-nu1}, $\GraphS{G}$ is $\F$-representable if and only if
$G$ has a nowhere-zero $(\F,+)$-flow. If $|\F|>2$, such a flow exists by
Lemma~\ref{lem:dirac-flows}. Over $\F_2$ it exists if and only if all degrees are even, by (F4). By
Remark~\ref{rem:nu1-examples}, $G$ is either $K_5$, which is $4$-regular, or $K_5$ minus an edge, a
wheel, or a graph obtained from $K_{3,p}$, and each of the latter has a vertex of degree~$3$.
\end{proof}

\begin{corollary}\label{cor:petersen}
The matroid $\GraphS{P_{10}}$ of the Petersen graph is representable over no field.
\end{corollary}
\begin{proof}
The Petersen graph is simple and $3$-connected, and it has two disjoint $5$-cycles, so
Theorem~\ref{thm:three-connected} applies. Explicitly, write its vertices as $u_i,v_i$, $i\in\mathbb{Z}_5$, with edges $u_iu_{i+1}$, $v_iv_{i+2}$ and $u_iv_i$. Contracting the matching
$\{u_0u_1,u_2u_3,v_0v_2,v_1v_4\}$ turns the two $5$-cycles into triangles, and the edges $u_0v_0$,
$u_3v_3$, $u_4v_4$ become a perfect matching between them; so Corollaries~\ref{cor:prism}
and~\ref{cor:loops}(3) exclude characteristic $2$.
\end{proof}

\section{Unicyclic graphs, the doubled triangle and Chow's remark}\label{sec:examples}

\begin{corollary}[Unicyclic graphs]\label{cor:unicyclic}
Let $G$ be a connected loopless graph with exactly one cycle $C$, of length $m\ge2$, let $D$ be its
reduced incidence matrix, and put
\[
 H_G=\left[\begin{array}{c|c} D & D\\ 0 & \one_C\end{array}\right]\in\F^{|V|\times2|E|}.
\]
\begin{enumerate}
\item If $\charF=2$, then $H_G$ represents $\GraphS{G}$.
\item If $m\ge3$, equivalently if $G$ is simple, then $\GraphS{G}$ is $\F$-representable if and only
if $\charF=2$.
\item If $m=2$, then $\GraphS{G}$ is representable over every field.
\end{enumerate}
\end{corollary}
\begin{proof}
(1) Here $\nuo(G)=1$, and the edges of $C$ are exactly the edges of $G$ that are not bridges. In
characteristic $2$, $\one_C$ is a flow on $G$, because every vertex of $C$ meets exactly two edges of
$C$. In Proposition~\ref{prop:lift} with $\lambda=\one_C$ we have $d_e=1$ for all $e$, so
$H_\lambda=H_G$. (2) and (3) Two parallel edges form a cycle, so $G$ is simple if and only if
$m\ge3$. If $m\ge3$, the edges of $C$ have no parallel mates, and Theorem~\ref{thm:oddchar}
excludes every characteristic $\neq2$. If $m=2$, the two edges of $C$ are parallel mates of each
other and are the only edges on cycles, so Theorem~\ref{thm:oddchar} gives representability in
every characteristic $\neq2$. Combine with (1).
\end{proof}

Since $|E|=|V|$ for a connected unicyclic graph, $\GraphS{G}$ is Lagrangian; over a finite field of characteristic $2$, $H_G$ is therefore the stabilizer matrix of a quantum stabilizer state. Adding a single edge to a unicyclic graph can destroy $\F_2$-representability: the doubled triangle below has $\nuo=1$, but its bridgeless part is not Eulerian.

\begin{example}\label{new:ex:unicyclic}
Let $n=|E|$, and label the cycle edges by $1,\ldots,m$ and the remaining edges by $m+1,\ldots,n$. For the unicyclic case in characteristic $2$, the matrix $M_0=[A\mid I_n]$, with $$A=\begin{bmatrix}
\mathcal{J}_m-I_m & \mathbf{0}_{m\times (n-m)}\\
\mathbf{0}_{(n-m)\times m} & I_{n-m}
\end{bmatrix},$$ where \(\mathcal{J}_m\) is the all-ones matrix and \(I_m\) is the identity matrix, represents $\mathcal{S}(G)$ when the cycle length $m$ is odd. When $m$ is even, interchanging the columns $1$ and $1^*$ gives a representation. 

\begin{proof}
The matrix $A$ is symmetric, so $M_0=[A\mid I_n]$ is isotropic and has rank $n$. Let $T$ be a transversal; put
\(Y=\{i\in[n]:i^*\in T\}\), \(U=[n]\setminus Y\), and \(
s=|U\cap[m]|\).
Eliminating the unit columns indexed by $Y$ gives
\[
\det M_0[T]
 =\pm\det A[U,U]
 =\pm\det(\mathcal{J}_s-I_s),
\]
where the determinant of the empty matrix is $1$. For $s\ge1$,
\[
\det(\mathcal{J}_s-I_s)=(-1)^{s-1}(s-1).
\]
Since $\operatorname{Char}\mathbb F=2$, this determinant is nonzero
exactly when $s$ is even. If $t=|Y\cap[m]|$ is the number of
starred cycle edges, then $s=m-t$. Hence, when $m$ is odd,
$\det M_0[T]\ne0$ exactly when $t$ is odd, as required by Chow's
rule.

When $m$ is even, swapping columns $1$ and $1^*$ replaces $Y$ by $Y\triangle\{1\}$, where \(\Delta\) is the symmetric difference, and therefore reverses this parity condition.
The swap preserves isotropy in characteristic two. Thus the resulting matrix has precisely the required admissible bases.
\end{proof}
\end{example}

\begin{corollary}\label{cor:doubled}
Let $\Kplus$ have vertices $1,2,3$, two parallel edges $e_1,e_2$ joining $1$ and $2$, and edges
$e_3=23$, $e_4=31$. Then $\GraphS{\Kplus}$ is $\F$-representable if and only if $\charF=2$ and
$|\F|>2$. In that case, for distinct $\alpha,\mu\in\F^\times$ and $\kappa=\alpha+\mu$, the matrix
\[ H_{\alpha,\mu}=\left[\begin{array}{cccc|cccc}
1&1&0&1&\alpha&\mu&0&\kappa\\
1&1&1&0&\alpha&\mu&\kappa&0\\
0&0&0&0&\alpha&\mu&\kappa&\kappa\end{array}\right], \]
with columns $e_1,\dots,e_4,e_1^*,\dots,e_4^*$ and the root $3$, represents $\GraphS{\Kplus}$.
Over $\F_4=\F_2(\vartheta)$, $\vartheta^2=\vartheta+1$, one may take $(\alpha,\mu)=(1,\vartheta)$ and
$\kappa=\vartheta^2$.
\end{corollary}
\begin{proof}
The edge $e_3$ lies on a cycle and has no parallel mate, so Theorem~\ref{thm:oddchar} excludes
characteristic $\neq2$. The vertices $1$ and $2$ have degree $3$ in $(\Kplus)^\circ=\Kplus$, so
Theorem~\ref{thm:flow} and (F4) exclude $\F_2$. If $\charF=2$ and $|\F|>2$, then
$\lambda=(\alpha,\mu,\kappa,\kappa)$ is a nowhere-zero flow: the sums at the vertices $1$, $2$ and
$3$ are $\alpha+\mu+\kappa=0$, $\alpha+\mu+\kappa=0$ and $\kappa+\kappa=0$, and $\kappa\neq0$ because
$\alpha\neq\mu$. Since $\nuo(\Kplus)=1$ and $\Kplus$ is bridgeless, Proposition~\ref{prop:lift}
applies with $d=\lambda$, and $H_\lambda=H_{\alpha,\mu}$, because
$\bd_{e_1}=\bd_{e_2}=\one_1+\one_2$, $\bd_{e_3}=\one_2$ and $\bd_{e_4}=\one_1$.
\end{proof}

\begin{remark}[Chow's remark]\label{rem:chow}
Chow \cite{Chow2003}, crediting a personal communication of N.~White, states that a graph with three
vertices and four edges, two of which form a double edge, has a non-representable graphical
symplectic matroid. A connected loopless graph with three vertices and four edges containing
exactly one pair of parallel edges is the doubled triangle $\Kplus$: the two remaining edges must
join the third vertex to both ends of the double edge. By Corollary~\ref{cor:doubled},
$\GraphS{\Kplus}$ is not representable over $\F_2$ or over any field of characteristic $\neq2$, but
it \emph{is} representable over $\F_4$, explicitly by $H_{1,\vartheta}$. Since the definition of
representability in \cite{BorovikGelfandWhite2003} allows every field, the remark holds precisely
with the qualification ``over $\F_2$ or over a field of characteristic $\neq2$''.
\end{remark}

\section{Rank enumerators: intrinsic versus enveloped}\label{sec:enumerators}

\subsection{The intrinsic corank--nullity polynomial}

\begin{definition}\label{def:TS}
Let $\cS$ be a symplectic matroid of rank $k$ on $J$. Its \emph{admissible corank--nullity
polynomial} and its \emph{restricted Tutte--Martin polynomial} are
\begin{equation}\label{eq:TS}
 T_{\cS}(u,v)=\sum_{A\in\Adm(J)}u^{\,k-r_{\cS}(A)}\,v^{\,|A|-r_{\cS}(A)},\qquad
 m_{\cS}(x)=\sum_{A\in\Adm_k(J)}(x-1)^{\,k-r_{\cS}(A)} .
\end{equation}
For $a\in J$ we also put
\[ T^{\circ}_{\cS,a}(u,v)=\sum_{A\in\Adm(J\setminus J_a)}u^{\,k-r_{\cS}(A)}\,v^{\,|A|-r_{\cS}(A)} . \]
\end{definition}

By Lemma~\ref{lem:rank-properties}(2) all exponents are nonnegative. Both polynomials depend only
on $\cS$. Since $|A|=k$ exactly when the two exponents in \eqref{eq:TS} coincide,
$m_{\cS}(x)=\sum_{j}c_{jj}(x-1)^{j}$, where $c_{ij}$ is the coefficient of $u^iv^j$ in $T_{\cS}$.
Moreover $T_{\cS}(1,1)=3^n$, $T_{\cS}(1,0)=|\cI(\cS)|$, $T_{\cS}(0,0)=|\cB(\cS)|$, and
$T_{\cS}(0,1)$ is the number of admissible sets containing a basis. For Lagrangian matroids
$m_{\cS}$ is Bouchet's restricted Tutte--Martin polynomial \cite{Bouchet2005}, and in general it is
the polynomial (13) of \cite{Sarvepalli2011}, which Sarvepalli defined for represented $\cS$ using
the ranks of column sets of a representation; the next proposition shows that the two agree.

\begin{proposition}\label{prop:rep-indep}
Let $H$ represent $\cS$, and let $M_H$ be the ordinary matroid on $J$ represented by the columns of $H$. Then $r_{M_H}(A)=r_{\cS}(A)$ for every $A\in\Adm(J)$. Consequently $\sum_{A\in\Adm_k(J)}(x-1)^{k-\rk H[A]}=m_{\cS}(x)$ and $\sum_{A\in\Adm(J)}u^{k-\rk H[A]}v^{|A|-\rk H[A]}=T_{\cS}(u,v)$; both are independent of the representation.
\end{proposition}
\begin{proof}
This is Lemma~\ref{lem:intrinsic}.
\end{proof}

\begin{proposition}\label{prop:contraction-rank}
Let $\cS$ be a symplectic matroid and $a\in J$. For every $X\in\Adm(J\setminus J_a)$,
\[ r_{\cS}(X\cup\{a\})=\begin{cases}1+r_{\cS/a}(X)&\text{if $a$ is not a loop of $\cS$,}\\
 r_{\cS}(X)&\text{if $a$ is a loop of $\cS$.}\end{cases} \]
\end{proposition}
\begin{proof}
If $a$ lies in no basis, then $|(X\cup\{a\})\cap B|=|X\cap B|$ for every basis $B$. Otherwise let
$T$ be a transversal containing $X\cup\{a\}$. By Lemmas~\ref{lem:restriction}
and~\ref{lem:rank-properties}(1), $r_{\cS}(X\cup\{a\})$ is the rank of $X\cup\{a\}$ in the matroid
$\cS|_T$, in which $\{a\}$ is independent. In an ordinary matroid,
$r(X\cup\{a\})=1+r_{/a}(X)$ for a non-loop $a$, where $r_{/a}$ is the rank function of the
contraction. The independent sets of $(\cS|_T)/a$ are the sets $I\subseteq T\setminus\{a\}$ with
$I\cup\{a\}\in\cI(\cS)$, that is, by Lemma~\ref{lem:minors}(1), the independent sets of $\cS/a$
contained in $T\setminus\{a\}$. Hence $r_{(\cS|_T)/a}(X)=r_{\cS/a}(X)$.
\end{proof}

\begin{theorem}\label{thm:threeterm}
Let $\cS$ be a symplectic matroid of rank $k$ on $J$, and let $a\in J$.
\begin{enumerate}
\item If neither $a$ nor $a^*$ is a loop, then
$T_{\cS}=T^{\circ}_{\cS,a}+T_{\cS/a}+T_{\cS/a^*}$.
\item If $a$ is a loop and $a^*$ is not, then $T_{\cS}=(1+v)\,T^{\circ}_{\cS,a}+T_{\cS/a^*}$.
\item If $a$ and $a^*$ are both loops, then $T_{\cS}=(1+2v)\,T^{\circ}_{\cS,a}$.
\end{enumerate}
\end{theorem}
\begin{proof}
Every admissible set either avoids $J_a$, or has the form $X\cup\{a\}$ or $X\cup\{a^*\}$ with
$X\in\Adm(J\setminus J_a)$, and these three cases are disjoint. The first case contributes
$T^{\circ}_{\cS,a}$. If $a$ is not a loop, then $k\ge1$, and Proposition~\ref{prop:contraction-rank}
gives $k-r_{\cS}(X\cup\{a\})=(k-1)-r_{\cS/a}(X)$ and
$|X\cup\{a\}|-r_{\cS}(X\cup\{a\})=|X|-r_{\cS/a}(X)$. Since $\cS/a$ is a symplectic matroid of rank
$k-1$ on $J\setminus J_a$ (Lemma~\ref{lem:minors}), the second case contributes $T_{\cS/a}$. If $a$ is a loop, the exponent
of $u$ is unchanged and that of $v$ increases by one, so the second case contributes
$v\,T^{\circ}_{\cS,a}$. The same applies to $a^*$.
\end{proof}

\begin{remark}\label{rem:deletion-fails}
For Lagrangian matroids and for graphical symplectic matroids, $T^{\circ}_{\cS,a}$ is again the
polynomial of a symplectic matroid on $J\setminus J_a$, up to a power of $u$
(Corollaries~\ref{cor:lagrangian}, \ref{cor:graph-recursion} and~\ref{cor:marked-class}). This fails
in general for $k<n$.
Let $\cS_0$ be the symplectic matroid of rank $2$ on $[3]\sqcup[3]^*$ represented over $\F_2$ by
\[ \left[\begin{array}{ccc|ccc}1&1&0&0&0&0\\0&0&1&1&1&0\end{array}\right], \]
with columns $1,2,3,1^*,2^*,3^*$; its bases are $\{1,3\},\{2,3\},\{1,2^*\},\{2,1^*\}$. The maximal
admissible independent subsets of $J\setminus J_1=\{2,3,2^*,3^*\}$ are $\{2,3\}$ and $\{2^*\}$, of
different sizes. Since all maximal independent sets of a symplectic matroid are bases and hence
have the same size, the restriction of $r_{\cS_0}$ to $\Adm(J\setminus J_1)$ is not the rank function
of a symplectic matroid on $J\setminus J_1$. Even the polynomial cannot be realized:
$T^{\circ}_{\cS_0,1}=u^2v+u^2+3uv+3u+1$. If this were $u^mT_{\cS'}$ for a symplectic matroid $\cS'$,
the constant term would force $m=0$, the degree in $u$ would force $r(\cS')=2$, and the constant
term would show that $\cS'$ has a single basis. Then $\cS'$ would have exactly two elements that are
not loops, whereas the coefficient of $u$, which counts the independent singletons, is $3$.
\end{remark}

\begin{corollary}\label{cor:lagrangian}
Let $\cS$ be a Lagrangian matroid on $J=[n]\sqcup[n]^*$ and $i\in[n]$. Let $\cS(i)$ be the
\emph{projection}, the Lagrangian matroid on $J\setminus J_i$ with bases $B\setminus J_i$,
$B\in\cB(\cS)$ \cite[Def.~5.1]{EFLS2024}.
\begin{enumerate}
\item If neither $i$ nor $i^*$ is a loop, then $T_{\cS}=T_{\cS/i}+T_{\cS/i^*}+u\,T_{\cS(i)}$.
\item If $i$ is a loop, then $T_{\cS}=(1+u+uv)\,T_{\cS/i^*}$; symmetrically if $i^*$ is a loop.
\item $T_{\cS}(u,v)=U_{\cS}(u,uv)$ and $m_{\cS}(x)=U_{\cS}(0,x-1)$, where $U_{\cS}$ is the
$U$-polynomial of the delta-matroid of $\cS$ \cite[Def.~1.4]{EFLS2024},
\cite[Def.~1.3]{Larson2025}. In particular $m_{\cS}(x)=\mathrm{Int}_{\cS}(x-1)$ for the interlace
polynomial $\mathrm{Int}_{\cS}(v)=U_{\cS}(0,v)$ of \cite{EFLS2024}.
\item If $\cS$ is represented over $\F_2$ by $[I_n\mid A_\Gamma]$, where $A_\Gamma$ is the adjacency
matrix of a simple graph $\Gamma$ on $[n]$, then $m_{\cS}(x)=q(\Gamma;x)$ is the interlace polynomial of
Arratia, Bollob\'as and Sorkin \cite{ArratiaBollobasSorkin2004}.
\end{enumerate}
\end{corollary}
\begin{proof}
Every basis of $\cS$ is a transversal, so it meets $J_i$ in exactly one element.
(1) For $A\in\Adm(J\setminus J_i)$, $|A\cap B|=|A\cap(B\setminus J_i)|$, so $r_{\cS}(A)=r_{\cS(i)}(A)$,
and $\cS(i)$ has rank $n-1$. Hence $T^{\circ}_{\cS,i}=u\,T_{\cS(i)}$, and
Theorem~\ref{thm:threeterm}(1) applies. (2) If $i$ is a loop, every basis contains $i^*$, so
$r_{\cS}(A)=r_{\cS/i^*}(A)$ for $A\in\Adm(J\setminus J_i)$, $T^{\circ}_{\cS,i}=u\,T_{\cS/i^*}$, and
Theorem~\ref{thm:threeterm}(2) gives $T_{\cS}=(1+v)u\,T_{\cS/i^*}+T_{\cS/i^*}$. (3) By
\cite[Def.~1.3]{Larson2025}, which agrees with \cite[Def.~1.4]{EFLS2024} by
\cite[Prop.~3.1]{Larson2025}, $U_{\cS}(u,v)=\sum_{A\in\Adm(J)}u^{n-|A|}v^{(|A|-g(A))/2}$, where
$g(A)=2r_{\cS}(A)-|A|$ by Remark~\ref{rem:rank}. Thus
$U_{\cS}(u,v)=\sum_Au^{n-|A|}v^{|A|-r_{\cS}(A)}$, and substituting $uv$ for $v$ gives $T_{\cS}$.
Setting $u=0$ keeps the terms with $|A|=n$, which gives $m_{\cS}(x)=U_{\cS}(0,x-1)$. (Items (1) and
(2), after this substitution, are the defining recursion \cite[Def.~1.4]{EFLS2024}.) (4) For a
transversal $T$ let $Y=\{j:j^*\in T\}$. The columns of $T$ are the unit vectors $\one_j$, $j\notin Y$,
and the columns of $A_\Gamma$ indexed by $Y$; eliminating with the unit vectors gives
$\rk H[T]=n-|Y|+\rk A_\Gamma[Y]$, where $A_\Gamma[Y]$ is the principal submatrix on $Y$. Hence
$m_{\cS}(x)=\sum_{Y\subseteq[n]}(x-1)^{|Y|-\rk A_\Gamma[Y]}$, which is $q(\Gamma;x)$ by the theorem of Aigner
and van der Holst \cite{AignerVanDerHolst2004}.
\end{proof}

\begin{corollary}[Deletion--contraction for graphical symplectic matroids]\label{cor:graph-recursion}
Let $G$ be a graph, $e$ an edge, and $\delta_e=r(\GraphS{G})-r(\GraphS{G\setminus e})$; thus
$\delta_e\in\{0,1\}$, and $\delta_e=1$ exactly when $G\setminus e$ has more acyclic components than
$G$.
\begin{enumerate}
\item If $e$ is not a loop, then
$T_{\GraphS{G}}=u^{\delta_e}\,T_{\GraphS{G\setminus e}}+2\,T_{\GraphS{G/e}}$.
\item If $e$ is a loop, then
$T_{\GraphS{G}}=(1+v)\,u^{\delta_e}\,T_{\GraphS{G\setminus e}}+T_{\GraphS{G}/e^*}$.
\end{enumerate}
\end{corollary}
\begin{proof}
Let $\cS=\GraphS{G}$ and $k=r(\cS)$. The value of $\delta_e$ follows from Lemma~\ref{lem:rankSG},
because deleting one edge creates at most one new acyclic component. Admissible subsets of
$J\setminus J_e$ are admissible sets for $G\setminus e$, and $G[A]=(G\setminus e)[A]$ for them;
hence $r_{\cS}(A)=r_{\GraphS{G\setminus e}}(A)$, and $T^{\circ}_{\cS,e}=u^{\delta_e}T_{\GraphS{G\setminus e}}$.
(1) Neither $e$ nor $e^*$ is a loop of $\cS$, since a single non-loop edge is a tree. By
Corollary~\ref{cor:loops}(2), $\cS/e=\GraphS{G/e}$. Let $s$ be an end of $e$. By
Lemma~\ref{lem:switching}, $\sigma_s$ maps $\cB(\cS)$ onto itself, and $\sigma_s(e)=e^*$; hence
$\cB(\cS/e^*)=\sigma_s(\cB(\cS/e))$. Since $\sigma_s$ permutes $\Adm(J\setminus J_e)$ preserving sizes,
$T_{\cS/e^*}=T_{\cS/e}$. Now apply Theorem~\ref{thm:threeterm}(1). (2) A loop without star is a cycle
with no starred edge, so $e$ is a loop of $\cS$, while $\{e^*\}$ is independent. Apply
Theorem~\ref{thm:threeterm}(2).
\end{proof}

In (2), $\GraphS{G}/e^*$ is in general not a graphical symplectic matroid: its independent sets are
the sets $X$ such that $X\cup\{e^*\}\in\cI(G)$, that is, the component of $G[X]$ containing the
vertex of the loop must be a tree. Matroids of this kind are the subject of the next subsection.

\subsection{Marked vertices}\label{subsec:marked}

\begin{definition}[Marked vertices]\label{def:marked}
Let $R\subseteq V$ be a set of \emph{marked} vertices of $G$. An admissible set $S$ is
\emph{independent in $\GraphS{G,R}$} if every component of $G[S]$ without a vertex of $R$ is good
(Definition~\ref{def:graphical}), and every component of $G[S]$ that meets $R$ is a tree containing
exactly one vertex of $R$. Let $\tau_R(G)$ be the number of components of $G$ that contain neither a
cycle nor a vertex of $R$.
\end{definition}

For $R=\emptyset$ these are the independent sets of $\GraphS{G}$, and $\tau_\emptyset(G)=\tau(G)$.

\begin{proposition}\label{prop:marked}
Let $G^R$ be obtained from $G$ by adding a loop $\ell_s$ at every $s\in R$, and let
$L_R=\{\ell_s^*:s\in R\}$. Then $L_R$ is independent in $\GraphS{G^R}$, and the independent sets of
$\GraphS{G^R}/L_R$ are the independent sets of Definition~\ref{def:marked}. Consequently these are
the independent sets of a symplectic matroid $\GraphS{G,R}$, and
\[ r(\GraphS{G,R})=|V|-|R|-\tau_R(G). \]
\end{proposition}
\begin{proof}
The nontrivial components of $G^R[L_R]$ are single loops with one starred edge, so $L_R$ is
independent. By Lemma~\ref{lem:minors}(1), applied $|R|$ times, an admissible set $S$ on the edges of
$G$ is independent in $\GraphS{G^R}/L_R$ if and only if $S\cup L_R$ is independent in $\GraphS{G^R}$.
Let $K$ be a component of $G[S]$ with $t$ marked vertices. In $G^R[S\cup L_R]$ it acquires $t$ starred
loops. If $t=0$, nothing changes. If $t\ge1$, the new component has $|E(K)|+t\ge|V(K)|-1+t$ edges, so
it is good if and only if it is unicyclic, that is, if and only if $K$ is a tree and $t=1$; its cycle
is then a starred loop. This proves the first two claims. By Lemmas~\ref{lem:minors}(1)
and~\ref{lem:rankSG}, $r(\GraphS{G,R})=r(\GraphS{G^R})-|R|=|V|-\tau(G^R)-|R|$, and the components of
$G^R$ without cycles are the components of $G$ without cycles and without marked vertices.
\end{proof}

On a transversal $T$, Proposition~\ref{prop:marked} exhibits $\GraphS{G,R}|_T$ as the contraction of
the frame matroid $\GraphS{G^R}|_{T\cup L_R}$ by negative loops, which behave like the half-edges of
\cite{Zaslavsky1991}. For an edge $e$ of $G$ put
\[ \delta^R_e=r(\GraphS{G,R})-r(\GraphS{G\setminus e,R})=\tau_R(G\setminus e)-\tau_R(G). \]
Deleting $e$ replaces the component $K$ of $G$ containing $e$ by one or two components. If $K$ has
no cycle and no marked vertex, neither have these, and if two of them have none, then neither has $K$.
Hence $\delta^R_e\in\{0,1\}$.

\begin{theorem}\label{thm:marked}
Let $R\subseteq V$, write $T_{G,R}=T_{\GraphS{G,R}}$, and let $e$ be an edge with ends $s$ and $t$,
where $s=t$ if $e$ is a loop.
\begin{enumerate}
\item If $s,t\in R$, then $T_{G,R}=(1+2v)\,T_{G\setminus e,R}$.
\item If $s\neq t$ and $\{s,t\}\not\subseteq R$, then
$T_{G,R}=u^{\delta^R_e}\,T_{G\setminus e,R}+2\,T_{G/e,R_e}$, where $R_e$ is the image of $R$ in $G/e$.
\item If $s=t\notin R$, then
$T_{G,R}=(1+v)\,u^{\delta^R_e}\,T_{G\setminus e,R}+T_{G\setminus e,R\cup\{s\}}$.
\end{enumerate}
If $G$ has no edges, then $T_{G,R}=1$. Thus (1)--(3) compute $T_{\GraphS{G}}=T_{G,\emptyset}$
recursively, removing one edge at a time.
\end{theorem}
\begin{proof}
An admissible set that avoids $J_e$ has the same independent subsets, hence the same rank, in
$\GraphS{G,R}$ and in $\GraphS{G\setminus e,R}$. Therefore
$T^{\circ}_{\GraphS{G,R},e}=u^{\delta^R_e}T_{G\setminus e,R}$.

(1) Here $e$ joins two marked vertices or is a loop at a marked vertex, so neither $\{e\}$ nor
$\{e^*\}$ is independent, and both are loops of $\GraphS{G,R}$. The components of $G$ and of
$G\setminus e$ that contain an end of $e$ contain a marked vertex, so $\delta^R_e=0$, and
Theorem~\ref{thm:threeterm}(3) gives (1).

(2) A single edge with at most one marked end is a tree, so neither $e$ nor $e^*$ is a loop. Let $X$
be admissible with $X\cap J_e=\emptyset$. Contracting $e$ maps the components of $G[X\cup\{e\}]$
bijectively onto those of $(G/e)[X]$. As in the proof of Lemma~\ref{lem:graph-contraction}, trees
correspond to trees and unicyclic components to unicyclic components with the same starred edges on
their cycles, and since $s$ and $t$ are not both marked, corresponding components contain the same
number of marked vertices. Hence $X\cup\{e\}$ is independent in $\GraphS{G,R}$ if and only if $X$ is
independent in $\GraphS{G/e,R_e}$, that is, $\GraphS{G,R}/e=\GraphS{G/e,R_e}$ by
Lemma~\ref{lem:minors}(1). The switching $\sigma_s$ of Lemma~\ref{lem:switching} preserves
independence in $\GraphS{G,R}$, because it does not change $G[S]$ and preserves the parity of the
number of starred edges on every cycle. As in the proof of Corollary~\ref{cor:graph-recursion},
$\sigma_s$ maps $\GraphS{G,R}/e$ onto $\GraphS{G,R}/e^*$, so both contractions have the polynomial
$T_{G/e,R_e}$, and Theorem~\ref{thm:threeterm}(1) gives (2).

(3) The unstarred loop $e$ is a cycle without starred edges in a component without marked vertices,
so $e$ is a loop of $\GraphS{G,R}$, whereas $\{e^*\}$ is independent. For $X$ as above,
$X\cup\{e^*\}$ is independent if and only if the component of $s$ in $G[X]$ is a tree without marked
vertices and the other components satisfy Definition~\ref{def:marked}, that is, if and only if $X$ is
independent in $\GraphS{G\setminus e,R\cup\{s\}}$. Hence $\GraphS{G,R}/e^*=\GraphS{G\setminus e,R\cup\{s\}}$,
and Theorem~\ref{thm:threeterm}(2) gives (3).

If $G$ has no edges, the empty set is the only admissible set and the rank is $0$.
\end{proof}

Call two symplectic matroids \emph{isomorphic} if a bijection between their ground sets that
commutes with the involution maps the bases of one onto the bases of the other.

\begin{corollary}\label{cor:marked-class}
Let $\mathcal M$ be the class of symplectic matroids isomorphic to some $\GraphS{G,R}$; it contains
all graphical symplectic matroids. Let $\cS\in\mathcal M$ be a symplectic matroid on $J$, and let
$a\in J$.
\begin{enumerate}
\item The independent sets of $\cS$ contained in $J\setminus J_a$ are those of a member
$\cS\setminus J_a$ of $\mathcal M$, and $T^{\circ}_{\cS,a}=u^{\,r(\cS)-r(\cS\setminus J_a)}\,T_{\cS\setminus J_a}$.
\item If $a$ is not a loop of $\cS$, then $\cS/a\in\mathcal M$.
\end{enumerate}
Hence Theorem~\ref{thm:threeterm} computes $T_{\cS}$ for every $\cS\in\mathcal M$ by a recursion
that stays inside $\mathcal M$.
\end{corollary}
\begin{proof}
We may assume that $\cS=\GraphS{G,R}$ and $a\in\{e,e^*\}$ for an edge $e$. Then
$\cS\setminus J_e=\GraphS{G\setminus e,R}$, and the formula for $T^{\circ}_{\cS,a}$ was shown at the
beginning of the proof of Theorem~\ref{thm:marked}. In case (1) of that theorem both elements of
$J_e$ are loops. In case (2), $\cS/e=\GraphS{G/e,R_e}$, and $\sigma_s$ is an isomorphism from
$\cS/e$ onto $\cS/e^*$. In case (3), $e$ is a loop and $\cS/e^*=\GraphS{G\setminus e,R\cup\{s\}}$.
\end{proof}

By Remark~\ref{rem:deletion-fails}, a class with property (1) cannot contain all symplectic
matroids.

\subsection{Envelopes}

\begin{definition}[Envelopes]\label{def:envelope}
Let $\cS$ be a symplectic matroid of rank $k$ on $J$. A \emph{pair-compatible envelope} of $\cS$ is
an ordinary matroid $M$ of rank $k$ on $J$ with $\cB(M)\cap\Adm_k(J)=\cB(\cS)$. These are the
enveloping matroids of \cite[Exercise~3.12.6]{BorovikGelfandWhite2003} (cf.\
\cite[Rem.~6.7]{EFLS2024}). An envelope $M$ is \emph{sheltering} if its admissible independent sets are exactly the independent sets of $\cS$ \cite[Sec.~4]{Bouchet1997}. The \emph{envelope polynomial} is
\[
 T_{\cS,M}(u,v)=\sum_{A\in\Adm(J)}u^{\,k-r_M(A)}v^{\,|A|-r_M(A)} .
\]
\end{definition}

\begin{proposition}\label{prop:envelope}
Let $M$ be a pair-compatible envelope of $\cS$. Then $r_M(A)\ge r_{\cS}(A)$ for every
$A\in\Adm(J)$, with equality for all $A$ if and only if $M$ is sheltering. Moreover
$T_{\cS,M}=T_{\cS}$ if and only if $M$ is sheltering.
\end{proposition}
\begin{proof}
Bases of $\cS$ are bases of $M$, so $A\cap B$ is $M$-independent for every $B\in\cB(\cS)$, which
gives the inequality. If equality holds everywhere and $A$ is admissible and $M$-independent, then
$r_{\cS}(A)=r_M(A)=|A|$, so $A\in\cI(\cS)$; conversely, if $M$ is sheltering, the $M$-independent
subsets of an admissible $A$ are the $\cS$-independent ones, so $r_M(A)=r_{\cS}(A)$. If
$r_M(A_0)>r_{\cS}(A_0)$ for some $A_0$, then, differentiating with respect to $u$,
\[\partial_uT_{\cS,M}(1,1)=\sum_A\big(k-r_M(A)\big)<\sum_A\big(k-r_{\cS}(A)\big)=\partial_uT_{\cS}(1,1). \qedhere\]
\end{proof}

\begin{corollary}\label{cor:sheltering-rep}
If $H$ represents $\cS$, then $M_H$ is a sheltering envelope of $\cS$ and $T_{\cS,M_H}=T_{\cS}$. In
particular, a symplectic matroid without a sheltering envelope is not representable over any
field.
\end{corollary}
\begin{proof}
$M_H$ has rank $k$, and its admissible bases are the admissible $k$-sets with nonzero determinant,
that is, the bases of $\cS$. It is sheltering by Lemma~\ref{lem:intrinsic}.
\end{proof}

\begin{example}[Duchamp's delta-matroid]\label{ex:duchamp}
Let $\cS$ be the Lagrangian matroid on $[4]\sqcup[4]^*$ with bases
\[\cB=\{1^*2^*3^*4^*,\ 1^*2^*3^*4,\ 1^*234^*,\ 12^*34^*,\ 123^*4^*,\ 1^*234,\ 12^*34,\ 123^*4,\ 1234\},\]
where $1^*2^*3^*4$ stands for $\{1^*,2^*,3^*,4\}$. Bouchet \cite[Sec.~4]{Bouchet1997} gives this
example, which he attributes to Duchamp, in his discussion of sheltering matroids; see also
\cite[Ex.~6.12]{EFLS2024}. A SAT enumeration (\ref{app:computations}) shows:
\begin{enumerate}
\item $\cS$ has exactly twelve pair-compatible envelopes; six of them have $49$ bases and six have
$53$ bases;
\item none of them is sheltering. For example, the rank-$4$ matroid $M_0$ on $J$ whose non-bases are
the $21$ four-element sets
\[\begin{array}{l}
1^*2^*33^*,\ 1^*2^*34,\ 1^*2^*34^*,\ 1^*22^*3,\ 1^*22^*3^*,\ 1^*23^*4,\ 1^*23^*4^*,\\
1^*233^*,\ 1^*244^*,\ 1^*3^*44^*,\ 11^*2^*3,\ 11^*23^*,\ 12^*3^*4,\ 12^*3^*4^*,\\
12^*44^*,\ 1234^*,\ 13^*44^*,\ 2^*3^*44^*,\ 22^*33^*,\ 23^*44^*,\ 33^*44^*
\end{array}\]
is an envelope with $49$ bases in which $\{1,2^*,3^*\}$ is independent, although no basis of $\cS$
contains it; so $r_{M_0}(\{1,2^*,3^*\})=3>2=r_{\cS}(\{1,2^*,3^*\})$;
\item the twelve envelopes have exactly two different polynomials. Let $M_1$ be the matroid whose
bases are those of $M_0$ together with the four inadmissible sets $11^*23^*$, $1^*233^*$,
$22^*33^*$ and $1^*22^*3^*$; it is one of the envelopes with $53$ bases. Then
\begin{align*}
T_{\cS}&=u^4+8u^3+3u^2v+24u^2+7uv+29u+9,\\
T_{\cS,M_0}&=u^4+8u^3+2u^2v+24u^2+7uv+30u+9,\\
T_{\cS,M_1}&=u^4+8u^3+u^2v+24u^2+7uv+31u+9,
\end{align*}
and nine of the envelopes, among them all six with $49$ bases, have the polynomial $T_{\cS,M_0}$,
while the other three have $T_{\cS,M_1}$.
\end{enumerate}
By (2) and Proposition~\ref{prop:envelope}, no envelope polynomial of $\cS$ equals $T_{\cS}$, and by
Corollary~\ref{cor:sheltering-rep}, $\cS$ is not representable over any field. By (3), the envelope
polynomial $T_{\cS,M}$ depends on the choice of $M$, whereas $T_{\cS}$ is an invariant of $\cS$. The
statement in \cite[Ex.~6.12]{EFLS2024} that no matroid on
$[4]\sqcup[4]^*$ has $\cB$ as its family of admissible bases is correct for sheltering matroids, the
notion relevant there (cf.\ \cite[Rem.~6.7]{EFLS2024}), but not for the enveloping matroids of
\cite[Exercise~3.12.6]{BorovikGelfandWhite2003}.
\end{example}

\begin{remark}\label{rem:uniform}
For $U^*_{k,n}$ one has $r(A)=\min(k,|A|)$ and there are $2^s\binom ns$ admissible $s$-sets, so
\[ T_{U^*_{k,n}}(u,v)=\sum_{s=0}^{k}2^s\binom ns u^{k-s}+\sum_{s=k+1}^{n}2^s\binom ns v^{s-k}. \]
Every envelope of $U^*_{k,n}$ is sheltering, because every admissible set of size at most $k$ is
independent in $U^*_{k,n}$; hence every envelope polynomial of $U^*_{k,n}$ equals
$T_{U^*_{k,n}}$. Let $P_{n,k}$ be the truncation to rank $k$ of the direct sum of the rank-one
uniform matroids on the pairs $J_1,\dots,J_n$; its independent sets are the admissible sets of size
at most $k$. It is an envelope of $U^*_{k,n}$, and it is the least envelope in the weak order: for
every envelope $M$, every independent set of $P_{n,k}$ lies in an admissible $k$-set, which is a
basis of $M$, so $r_{P_{n,k}}\le r_M$. For $k\ge2$ its rank-one flats are the inadmissible pairs
$J_i$ (and for $k=1$ its only rank-one flat is $J$), so $P_{n,k}$ is not an admissible matroid in
the sense of Raz, whose atoms are admissible \cite[Lemma~2.3]{Raz2024}. Relaxing an inadmissible
circuit-hyperplane of an envelope \cite[Prop.~1.5.14]{Oxley2011} produces another envelope with the
same ranks on admissible sets.
\end{remark}

\section{Open problems}\label{sec:open}

The normal form (Theorem~\ref{thm:normalform}) reduces the representation problem for
$\GraphS{G}$ to the choice of a flow, linear conditions on $Z_1$, and the prescribed vanishing and
nonvanishing of maximal minors (Remark~\ref{rem:decision}). The problem is solved in characteristic
$\neq2$, in characteristic $2$ for loopless graphs without two vertex-disjoint cycles, and over every
field for $3$-connected simple graphs (Theorem~\ref{thm:three-connected}). Characteristic $2$ in
general remains open.

\begin{question}\label{q:char2}
Characterize the graphs $G$ for which $\GraphS{G}$ is representable over a field of characteristic
$2$ with at least $8$ elements, where flows impose no restriction. If $\GraphS{G}$ is not
representable over such a field, does $G$ have a forest $F$ such that $G/F$ contains a rigidly linked
pair of cycles, as in Proposition~\ref{prop:prism-subdivision}?
\end{question}

\begin{question}\label{q:sums}
Over fields of characteristic $2$ with at least $4$ elements, is representability of $\GraphS{G}$
preserved under $1$-sums and under joining two graphs by a bridge? The graphs $X_1$ and $X_5$ of
Example~\ref{ex:char2} are instances of a positive answer, whereas $X_5$ shows that $1$-sums do not
preserve representability over $\F_2$.
\end{question}

\begin{question}\label{q:complexity}
By Theorem~\ref{thm:char2-nu1}, $\F_4$-representability of $\GraphS{G}$ is equivalent to the
existence of a nowhere-zero $4$-flow on $G^\circ$ when $\nuo(G)\le1$. Is deciding
$\F_4$-representability of $\GraphS{G}$ NP-hard in general, as deciding $3$-edge-colourability of
cubic graphs is \cite{Holyer1981}?
\end{question}

\appendix

\section{Computations}\label{app:computations}

Only Examples~\ref{ex:char2}(2),(3),(5) and~\ref{ex:duchamp} depend on computer calculations; all other statements are proved in the text.

\paragraph{Representability searches} For a connected graph $G$ with a fixed orientation and root
and a finite field $\F_q$, the program enumerates all $\F_q$-flows $z$ that are nonzero exactly on
the non-bridges, up to a scalar factor (the last row of a normal form may be rescaled). For each
such $z$ it solves the linear system in the entries of $Z_1$ formed by (N3) and by the relations
obtained, as in the proof of Lemma~\ref{lem:twostar}, from every dependent admissible set with at
most two starred elements whose unstarred part is a forest. By Theorem~\ref{thm:normalform}, every
representation in normal form arises in this way. If the solution space has at most $2^{18}$
elements, all solutions are enumerated; otherwise random solutions are sampled, which can only
establish representability. Each candidate $H$ is compared with Chow's rule on all admissible
$|V|$-sets. The results are as follows.
\begin{itemize}
\item $X_5$ over $\F_2$: the only nowhere-zero $\F_2$-flow is $\one_E$; the linear system has a
solution space of dimension $16$, and none of its $2^{16}$ elements gives a representation. Hence
$\GraphS{X_5}$ is not $\F_2$-representable.
\item $X_4$ over $\F_4$: there are $24$ flows up to scalars, each with a solution space of dimension
$5$, and none of the $24\cdot4^5$ candidates is a representation, in accordance with
Theorem~\ref{thm:rigid}.
\item $X_2$ over $\F_4$ and $\F_8$, $X_3$ over $\F_4$ and $X_5$ over $\F_4$: the representations
below were found and then checked against Chow's rule on all admissible $|V|$-sets.
\end{itemize}
We write $\F_4=\F_2(\vartheta)$ with $\vartheta^2=\vartheta+1$ and $\bar\vartheta=\vartheta^2$, and
$\F_8=\F_2(\eta)$ with $\eta^3=\eta+1$; in the matrix over $\F_8$, an entry $c=c_0+2c_1+4c_2$
($c_j\in\{0,1\}$) stands for $c_0+c_1\eta+c_2\eta^2$. In each matrix the columns are indexed by the
edges in the order listed in Example~\ref{ex:char2}, first unstarred and then starred, the root is
the vertex $0$, and the first $|V|-1$ rows are indexed by the vertices $1,\dots,|V|-1$.
{\footnotesize
\[
X_2,\ \F_4:\ \left[\begin{array}{cccccccc|cccccccc}
1&1&0&0&0&0&0&0&0&0&1&\bar\vartheta&\bar\vartheta&1&\vartheta&\vartheta\\
0&1&1&0&0&0&0&1&\bar\vartheta&0&0&\vartheta&\vartheta&\bar\vartheta&1&0\\
0&0&0&1&0&1&1&0&\bar\vartheta&\bar\vartheta&1&\vartheta&\bar\vartheta&\bar\vartheta&1&\vartheta\\
0&0&0&1&1&0&0&0&0&0&0&1&1&0&0&0\\
0&0&0&0&1&1&0&1&1&1&\vartheta&\vartheta&\bar\vartheta&1&1&\vartheta\\
0&0&0&0&0&0&0&0&1&1&\vartheta&1&1&\vartheta&\bar\vartheta&\bar\vartheta
\end{array}\right]
\]
\[
X_2,\ \F_8:\ \left[\begin{array}{cccccccc|cccccccc}
1&1&0&0&0&0&0&0&4&4&1&5&5&1&4&4\\
0&1&1&0&0&0&0&1&6&7&5&1&1&2&3&2\\
0&0&0&1&0&1&1&0&5&5&1&2&3&4&6&4\\
0&0&0&1&1&0&0&0&6&6&7&7&7&7&1&1\\
0&0&0&0&1&1&0&1&7&7&5&0&1&2&1&3\\
0&0&0&0&0&0&0&0&1&1&2&1&1&2&3&3
\end{array}\right]
\]
\[
X_3,\ \F_4:\ \left[\begin{array}{cccccc|cccccc}
1&1&0&0&0&1&0&0&1&\vartheta&1&0\\
0&0&1&1&1&0&0&0&\vartheta&\bar\vartheta&1&\vartheta\\
0&0&1&1&0&1&0&0&1&\vartheta&0&\bar\vartheta\\
0&0&0&0&0&0&1&\vartheta&1&\vartheta&\bar\vartheta&\bar\vartheta
\end{array}\right]
\]
\[
X_5,\ \F_4:\ \left[\begin{array}{ccccccccc|ccccccccc}
1&1&0&0&0&0&0&0&0&0&0&\vartheta&\vartheta&\vartheta&\vartheta&\vartheta&\vartheta&\vartheta\\
0&1&1&1&0&1&0&0&0&1&\bar\vartheta&\bar\vartheta&\bar\vartheta&\vartheta&\bar\vartheta&0&0&0\\
0&0&0&1&1&0&0&0&0&0&0&0&1&1&0&0&0&0\\
0&0&0&0&1&1&1&0&1&0&0&0&\bar\vartheta&\vartheta&\vartheta&\bar\vartheta&1&\bar\vartheta\\
0&0&0&0&0&0&1&1&0&0&0&0&0&0&0&\vartheta&\vartheta&0\\
0&0&0&0&0&0&0&1&1&0&0&0&0&0&0&0&\vartheta&\vartheta\\
0&0&0&0&0&0&0&0&0&1&1&1&1&1&1&1&1&1
\end{array}\right]
\]}

\paragraph{Envelopes of Duchamp's matroid} Rank-$4$ matroids on the $8$-element set $J$ are
encoded by $\binom84=70$ Boolean variables, one for each $4$-subset, with unit clauses fixing the
$16$ admissible $4$-sets according to $\cB$ and clauses expressing the basis exchange axiom. All
solutions were enumerated with the SAT solver Glucose by adding a blocking clause after each
solution, and each solution was checked independently. There are $12$ solutions, six with $49$ and
six with $53$ bases, and none is sheltering. Adding the clauses that make every admissible set not
contained in a member of $\cB$ dependent gives an unsatisfiable instance. A sheltering matroid of
rank larger than $4$ would truncate to a sheltering envelope of rank $4$, so $\cS$ has no sheltering
matroid of any rank. The polynomials $T_{\cS,M}$ of the twelve envelopes were computed from their
bases; they take the two values listed in Example~\ref{ex:duchamp}(3).

\paragraph{Further checks} The computer programs are also tested, on random and on structured
instances, the general statements of the paper: Lemmas~\ref{lem:restriction}, \ref{lem:minors},
\ref{lem:signed} (with loops) and~\ref{lem:intrinsic}; the constructions of
Theorem~\ref{thm:oddchar} (with loops), Proposition~\ref{prop:lift} and
Corollaries~\ref{cor:unicyclic} and~\ref{cor:doubled}; the flow facts used in
Corollary~\ref{cor:flow-small}; the cores in Example~\ref{ex:char2} and the matrix of its item~(1),
on all $448$ admissible $6$-sets; the forest of Proposition~\ref{prop:prism-subdivision} on random
prism subdivisions, the flows of Lemma~\ref{lem:dirac-flows}, the list of
Remark~\ref{rem:nu1-examples} against all $3$-connected simple graphs with at most seven vertices,
and Corollary~\ref{cor:petersen}; Theorem~\ref{thm:threeterm}, Corollaries~\ref{cor:lagrangian}
and~\ref{cor:graph-recursion}, and Proposition~\ref{prop:marked} and Theorem~\ref{thm:marked} on
random graphs with marked vertices. An exhaustive search also confirms that the prism is not
$\F_4$-representable, as Corollary~\ref{cor:prism} predicts. These checks are not used in the
proofs.

\bibliographystyle{elsarticle-num}
\bibliography{references}

\end{document}